\documentclass[11pt,a4paper]{article}

\usepackage{cite}
\usepackage[cmex10]{amsmath}
\usepackage{amsmath}
\usepackage{amssymb}
\usepackage{booktabs}
\usepackage{color}
\usepackage{times}
\usepackage{paralist}
\usepackage{subfigure}
\usepackage{float}
\usepackage{epstopdf}
\usepackage{authblk}
\usepackage{graphicx}

\usepackage[utf8]{inputenc}
\usepackage[cmex10]{amsmath}
\usepackage{amssymb}
\usepackage{mathtools}
\usepackage{subfigure}
\usepackage{float}
\usepackage{booktabs}
\usepackage{times}
\usepackage{paralist}
\usepackage{cite}
\usepackage{algorithm}
\usepackage{algpseudocode}
\usepackage{epstopdf}
\usepackage{soul} 
\usepackage{color}
\usepackage{eurosym}
\usepackage[printonlyused]{acronym}
\graphicspath{
    {FiguresBldsExp/}
}
\DeclareGraphicsExtensions{.pdf,.png,.jpg,.jpeg,.eps}

\usepackage{amsthm}

\definecolor{Red}{rgb}{1,0,0}
\definecolor{Green}{rgb}{0,1,0}
\definecolor{Blue}{rgb}{0,0,1}

\begin{document}
\title{Experimental Demonstration of Frequency Regulation by Commercial Buildings -- Part II: Results and Performance Evaluation}
\author[1]{Evangelos~Vrettos\thanks{vrettos@eeh.ee.ethz.ch}}
\author[2]{Emre~C.~Kara\thanks{eckara@slac.stanford.edu}}
\author[3]{Jason~MacDonald\thanks{jsmacdonald@lbl.gov}}
\author[1]{G\"{o}ran~Andersson\thanks{andersson@eeh.ee.ethz.ch}}
\author[4]{Duncan~S.~Callaway\thanks{dcal@berkeley.edu}}
\affil[1]{Power Systems Laboratory, ETH Zurich, Switzerland}
\affil[2]{National Accelerator Laboratory (SLAC), California, US}
\affil[3]{Grid Integration Group, Lawrence Berkeley National Laboratory (LBNL), California, US}
\affil[4]{Energy and Resources Group, University of California, Berkeley, US}
\maketitle

\begin{abstract}
This paper is the second part of a two-part series presenting the results from an experimental demonstration of frequency regulation in a commercial building test facility. In Part I, we developed relevant building models and designed a hierarchical controller for reserve scheduling, building climate control and frequency regulation.

In Part II, we introduce the communication architecture and experiment settings, and present extensive experimental results under frequency regulation. More specifically, we compute the day-ahead reserve capacity of the test facility under different assumptions and conditions. Furthermore, we demonstrate the ability of model predictive control to satisfy comfort constraints under frequency regulation, and show that fan speed control can track the fast-moving RegD signal of the Pennsylvania, Jersey, and Maryland Power Market (PJM) very accurately. In addition, we report the observed effects of frequency regulation on building control and provide suggestions for real-world implementation projects. Our results show that hierarchical control is appropriate for frequency regulation from commercial buildings.
\end{abstract}

\newpage
\section*{Acronyms}

\acrodefplural{RES}[RES's]{Renewable Energy Sources}
\acrodefplural{AS}[AS's]{Ancillary Services}

\begin{acronym}[swissgrid]
\acro{ACE}{Area Control Error}
\acro{AGC}{Automatic Generation Control}
\acro{AHU}{Air Handling Unit}
\acro{AMI}{Advanced Metering Infrastructure}
\acro{AMR}{Advanced Metering Reading}
\acro{AS}{Ancillary Service}
\acrodefplural{AS}[AS]{Ancillary Services}
\acro{AVR}{Automatic Voltage Regulator}
\acro{BAS}{Building Automation System}
\acro{BESS}{Battery Energy Storage System}
\acro{BE}{Balancing Energy}
\acro{BG}{Balance Group}
\acro{CAISO}{California Independent System Operator}
\acro{CDF}{Cumulative Distribution Function}
\acro{COP}{Coefficient of Performance}
\acro{CWS}{Central Working Station}
\acro{DA}{Asymmetric Reserves at a Building with Daily Duration}
\acro{DAE}{Differential-Algebraic Equation}
\acro{DN}{Distribution Network}
\acro{DR}{Demand Response}
\acro{DSA}{Symmetric Reserves in Aggregate with Daily Duration}
\acro{DSB}{Symmetric Reserves at a Building with Daily Duration}
\acro{DSM}{Demand Side Management}
\acro{DSO}{Distribution System Operator}
\acro{EKZ}{the utility of the Kanton Zurich}
\acro{ENTSO-E}{European Network of Transmission System Operators for Electricity}
\acro{ERCOT}{Electric Reliability Council of Texas}
\acro{EWH}{Electric Water Heater}
\acro{FIT}{Feed-in Tariff}
\acro{FLEXLAB}{Facility for Low Energy eXperiments}
\acro{GPU}{Graphics Processing Unit}
\acro{HA}{Asymmetric Reserves at a Building with Hourly Duration}
\acro{HP}{Heat Pump}
\acro{HSA}{Symmetric Reserves in Aggregate with Hourly Duration}
\acro{HSB}{Symmetric Reserves at a Building with Hourly Duration}
\acro{HVAC}{Heating, Ventilation and Air-Conditioning}
\acro{HVDC}{High-Voltage Direct Current}
\acro{HV}{High-Voltage}
\acro{IRA}{Integrated Room Automation}
\acro{KiBaM}{Kinetic Battery Model}
\acro{LBNL}{Lawrence Berkeley National Laboratory}
\acro{LP}{Linear Program}
\acro{LV}{Low-Voltage}
\acro{MAE}{Mean Absolute Error}
\acro{MAPE}{Mean Absolute Percentage Error}
\acro{MHSE}{Moving Horizon State Estimation}
\acro{MILP}{Mixed-Integer Linear Program}
\acro{MLD}{Mixed Logical Dynamical}
\acro{MPC}{Model Predictive Control}
\acro{MPPT}{Maximum Power Point Tracking}
\acro{MV}{Medium-Voltage}
\acro{NPV}{Net Present Value}
\acro{ODE}{Ordinary Differential Equation}
\acro{OLTC}{On-Load Tap-Changer}
\acro{OPF}{Optimal Power Flow}
\acro{PC}{Power Constraints}
\acro{PCC}{Point of Common Coupling}
\acro{PDE}{Partial Differential Equation}
\acro{PDF}{Probability Distribution Function}
\acro{PEC}{Power and Energy Constraints}
\acro{PEV}{Plug-in Electric Vehicle}
\acro{PF}{Power Flow}
\acro{PFC}{Primary Frequency Control}
\acro{PI}{Proportional-Integral}
\acro{PID}{Proportional-Integral-Derivative}
\acro{PJM}{Pennsylvania, Jersey, and Maryland Power Market}
\acro{PLC}{Power Line Communication}
\acro{PSS}{Power System Stabilizer}
\acro{PV}{Photovoltaics}
\acro{QP}{Quadratic Program}
\acro{RBC}{Rule-based Control}
\acro{RC}{Resistance-capacitance}
\acro{REF}{Refrigerator}
\acro{RES}{Renewable Energy Source}
\acrodefplural{RES}[RES]{Renewable Energy Sources}
\acro{RHC}{Receding Horizon Control}
\acro{RMSE}{Root Mean Squared Error}
\acro{RPC}{Robust Problem with Power Constraints}
\acro{RPEC}{Robust Problem with Power and Energy Constraints}
\acro{SAT}{Supply Air Temperature}
\acro{SC}{Slab Cooling}
\acro{SFC}{Secondary Frequency Control}
\acro{SLP}{Sequential Linear Programming}
\acro{SM}{Smart Meter}
\acro{SOC}[SoC]{State of Charge}
\acro{SOH}[SoH]{State of Health}
\acro{SPC}{Stochastic Problem with Power Constraints}
\acro{SPEC}{Stochastic Problem with power and Energy Constraints}
\acro{TABS}{Thermally Activated Building Systems}
\acro{TCL}{Thermostatically Controlled Load}
\acro{TSO}{Transmission System Operator}
\acro{sLP}[SLP]{Stochastic Linear Program}
\acro{swissgrid}{Swiss Transmission System Operator}
\acro{UL}{Uncontrollable Load}
\acro{VAV}{Variable Air Volume}
\acro{VFD}{Variable Frequency Drive}
\end{acronym}

\section{Introduction}
In Part I of this two-part paper, we performed a detailed literature review on theoretical, simulation-based, and experimental work on frequency regulation with commercial buildings. Furthermore, we presented the test facility for our experiment (FLEXLAB), developed relevant building models, and designed a hierarchical controller for reserve scheduling (level $1$), building climate \ac{MPC} (level $2$) and frequency regulation (level $3$).

In Section~\ref{exp_preparation} of Part II, we summarize the control and communication architecture, as well as the experiment settings. In Sections~\ref{results_lv1}, \ref{results_lv2}, and \ref{results_lv3} we report extensive experimental results in FLEXLAB for each level of control over a period of one week. We summarize some important findings and suggestions for future work in Section~\ref{sec_learnings}, whereas Section~\ref{sec_conclusion} concludes. 
\section{Preparation of the Experiment} \label{exp_preparation}
\subsection{Communication Architecture}

\begin{figure}[t]
\centering \includegraphics[width=1.05\textwidth, height=1.3in]{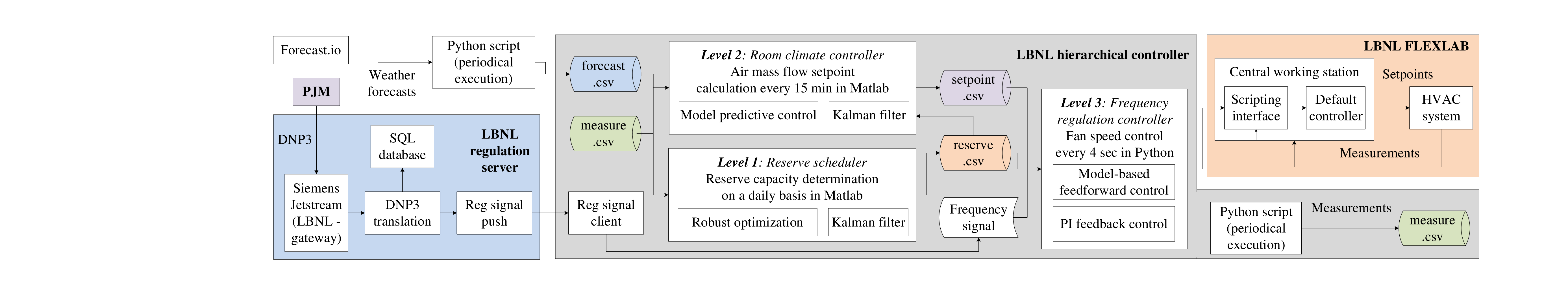}
\caption{The developed control and communication architecture for building climate control and frequency regulation in the FLEXLAB test facility.} \label{fig:comm_arch}
\end{figure}

We implement the reserve scheduling problem of level $1$ (solved once a day) and the MPC of level $2$ (solved every $15$ minutes) in Matlab, whereas we calculate the fan speed setpoints of level $3$ in Python (every $4$ seconds) and communicate them to the \ac{CWS} of FLEXLAB. We used a file-based communication between Python and Matlab based on comma-separated-values (csv) files.

The reserve scheduler stores the computed reserve capacity in the reserve.csv file. A Python script periodically queries the CWS and stores the building measurements in the measure.csv file. Another Python script periodically queries the publicly available database of forecast.io and stores the weather forecasts in the forecast.csv file.\footnote{Only ambient temperature forecasts are obtained from forecast.io. The solar radiation forecasts are obtained from a clear-sky radiation model, which turned out to be sufficient for the weather conditions during the experiment.} The MPC's feedback from the building is obtained from measure.csv and the weather forecasts from forecast.csv. The optimal air flow rate setpoint calculated by Matlab is stored in the setpoint.csv file. The fan speed setpoint is determined in Python by accessing the setpoint.csv and reserve.csv files, and based on the frequency regulation signal.

Most of the experiment was performed using archived data of the RegD signal from the \ac{PJM} from December $2012$ to January $2013$. Although the signal was available with a resolution of $2$~seconds, we down-sampled it to $4$~seconds due to the expected communication delays. In addition, a connection with PJM was established based on the DNP3 protocol and using a Siemens Jetstream gateway that provided us with the RegD signal in real-time. At the FLEXLAB side, the received data were translated, saved in an SQL database, and pushed by a ``RegD signal server'' to a ``RegD signal client''. The complete communication architecture from PJM to FLEXLAB is graphically shown in Fig.~\ref{fig:comm_arch}. However, network issues at FLEXLAB made the connection unreliable, and therefore we chose to run the live connection with PJM only for one continuous hour.

\subsection{Experiment Settings} \label{exp_settings}
\begin{figure}[t]
\centering \includegraphics[width=0.95\textwidth]{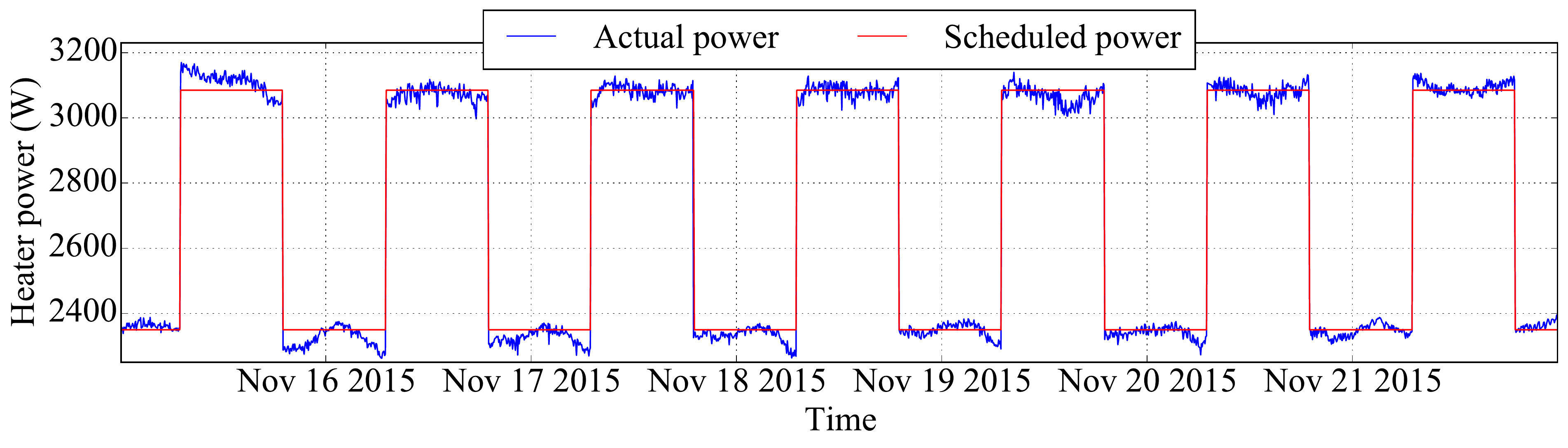}
\caption{The heater schedule and the actual power consumption. The heat gain is high during working hours and low during non-working hours.}
\label{fig:heater_plan}
\end{figure}

Since FLEXLAB is not occupied, we emulated the internal heat gains from occupants and equipment using electric heaters as plug loads. The total internal heat gain in both cells was kept lower than the chiller cooling capacity. The heaters' consumption profile was fixed according to the red curve of Fig.~\ref{fig:heater_plan} using digital timer sockets. The actual heater power (blue curve) fluctuates around the profile due to voltage variations.

Before the start of the experiment, we fixed the manually controlled inlet dampers in the rooms to fully open positions. In addition, we fixed the return air damper to a $100\%$ opening and the outside air damper to a $0\%$ opening, i.e., the return air was fully recirculated. The speeds of primary and secondary chilled water pumps were fixed to $75\%$ and $100\%$ of their rated speeds, respectively. Moreover, we deactivated the existing floor heating system and the heating coil at the \ac{AHU}.

\begin{figure}[t]
\centering \includegraphics[width=0.95\textwidth]{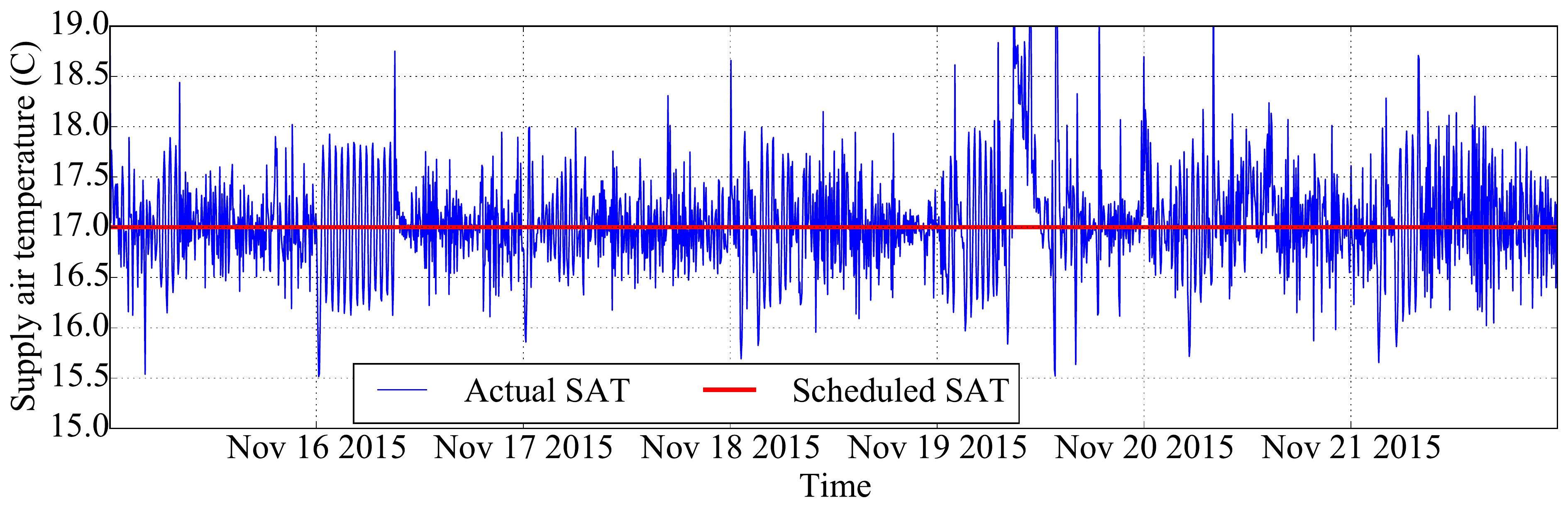}
\caption{The SAT setpoint and actual values during one week.}
\label{fig:sat_sch}
\end{figure}

We set the temperature comfort zone to $21-25^\circ$C during working hours. An existing \ac{PI} controller regulates the \ac{SAT} to $17^\circ$C by controlling the position of a cooling valve. The gains of this controller had been tuned for a conventional building operation; therefore, we modified them to achieve a tighter control and reduce the fluctuations of SAT around its setpoint during frequency regulation. The resulting SAT profile during the experiment is shown in Fig.~\ref{fig:sat_sch}. The mean deviation from the SAT setpoint is $0.05^\circ$C (there is a small bias to larger SAT values) and the mean absolute deviation is $0.37^\circ$C.

Recall that the facility has two building cells with identical construction: cell 1A is used for the frequency regulation experiment, whereas cell 1B serves as a benchmark. Applying the same air flow rate in both cells and recording the temperature, we verified that the two cells are thermally very similar. However, we observed that that the same fan speed setpoint induces a slightly different air flow rate in the two cells due to small differences in the \acp{AHU}. To compensate for this, we fitted different fan models for the two cells (the parameters for cell 1A are given in \cite[Table IV]{VrettosTSG2016Exp_p1}).

The electricity cost was assumed equal to $c_k = 0.18$~\euro/kWh, whereas the reserve capacity payment was fixed to a $10\%$ higher value, i.e., $\lambda_k=0.198$~\euro/kWh. The goal of this experiment is to demonstrate the technical feasibility of reserve provision from commercial buildings; therefore, we chose a relatively high capacity payment to incentivize reserve provision.

\subsection{Experiment Plan} \label{experiment_plan}
The experiment was organized into two parts. The first part took place from $15$ to $18$ November $2015$ and relied on an ``older'' building model identified with data from June-July $2015$ (see \cite[Table II]{VrettosTSG2016Exp_p1}). On $19$ November the experiment was paused and a new building model was identified using the recently collected data (see \cite[Table III]{VrettosTSG2016Exp_p1}), which was used in the second part of the experiment from $20$ to $21$ November.

\section{Reserve Scheduling (Level 1)} \label{results_lv1}
\begin{figure}[t]
\centering
\begin{minipage}{0.49\linewidth}%
\centering \includegraphics[width=1\textwidth]{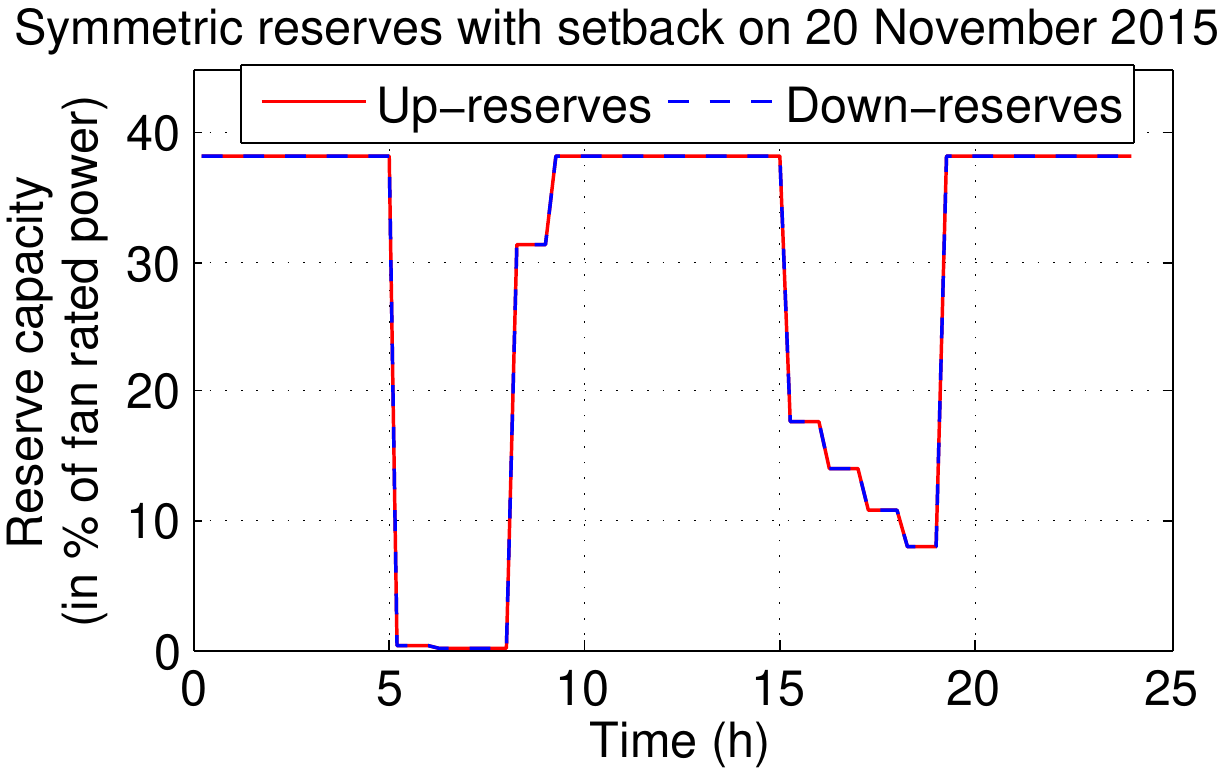}
\end{minipage}
\begin{minipage}{0.49\linewidth}
\centering \includegraphics[width=1\textwidth]{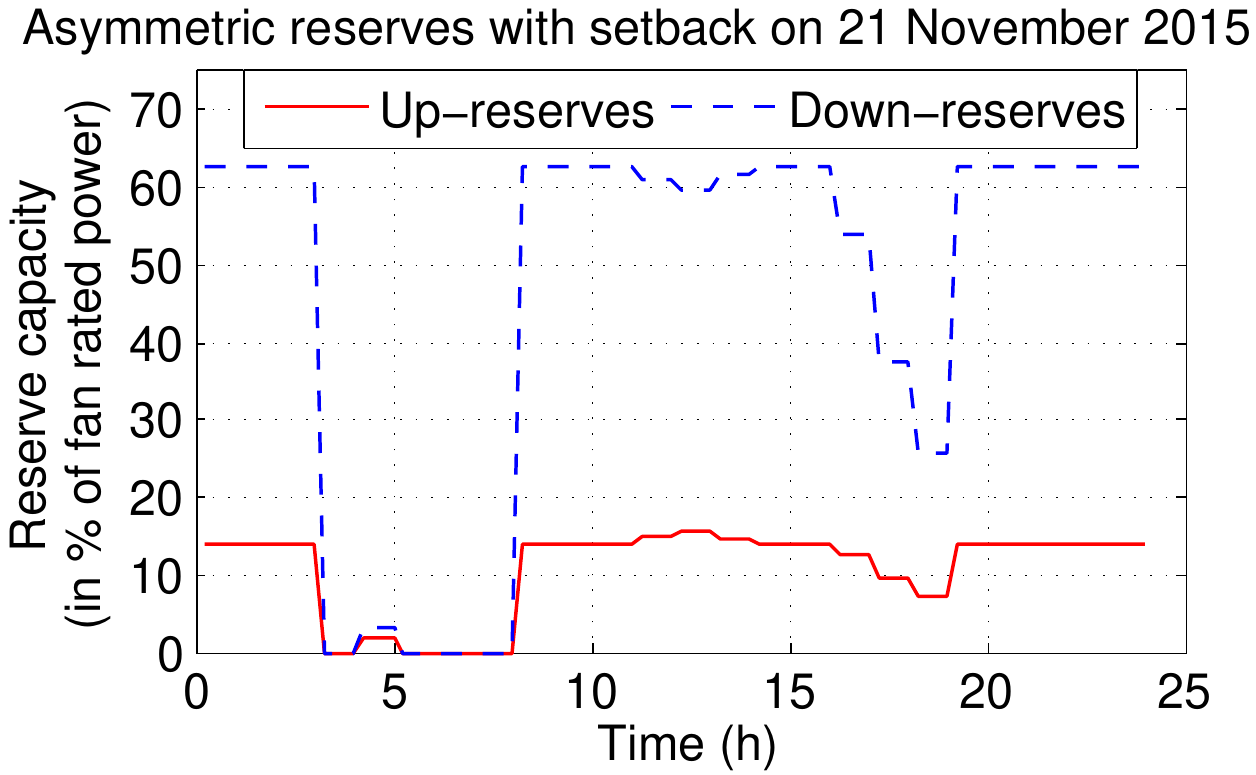}
\end{minipage}
\caption{The hourly reserve capacities as a percentage of the nominal fan power for $20$ November $2015$ (left) and $21$ November $2015$ (right).}
\label{fig:res_cap_sym}
\end{figure}

\begin{figure}[t]
\centering
\begin{minipage}{0.49\linewidth}%
\centering \includegraphics[width=1\textwidth]{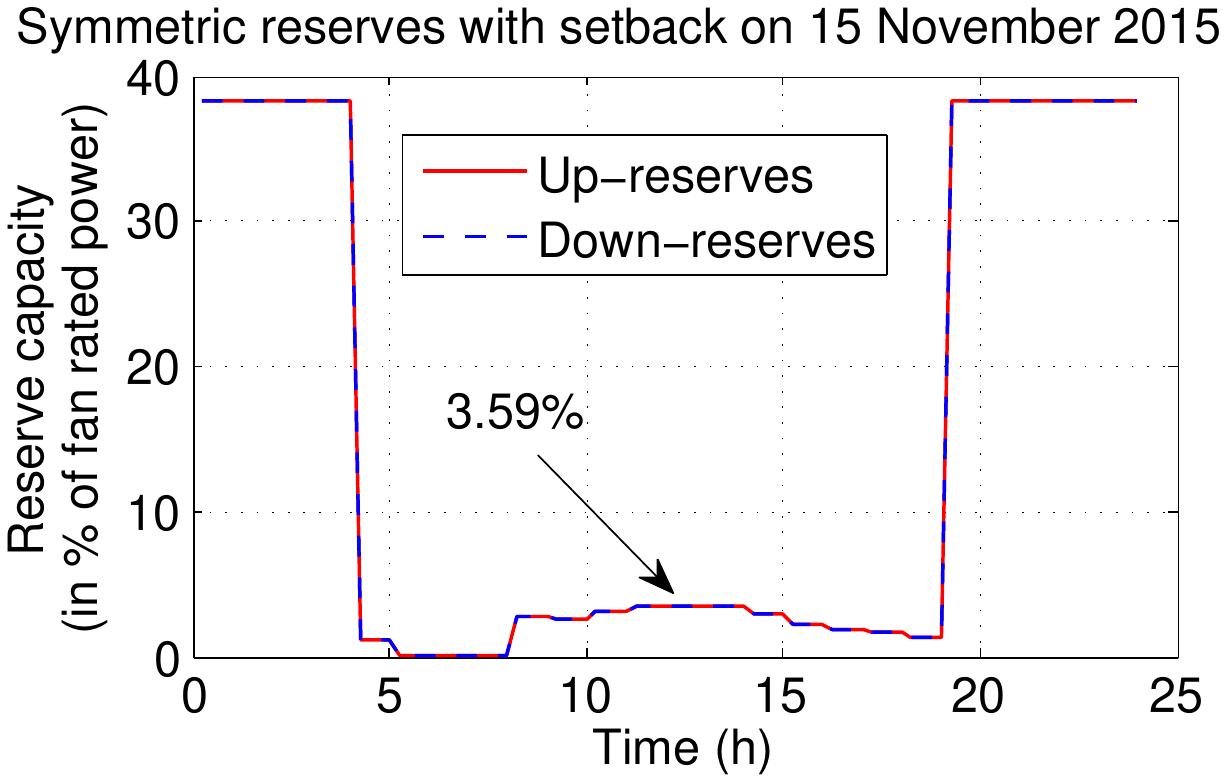}
\end{minipage}
\begin{minipage}{0.49\linewidth}
\centering \includegraphics[width=1\textwidth]{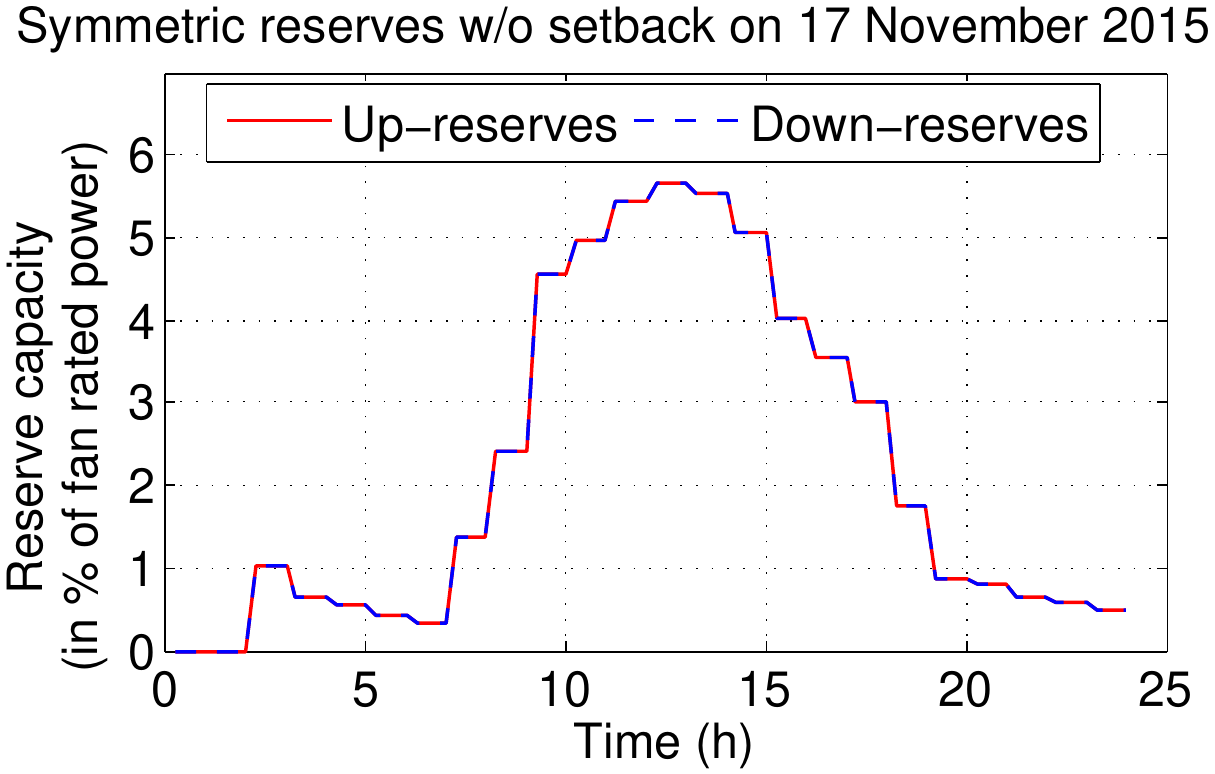}
\end{minipage}
\caption{The hourly reserve capacities as a percentage of the nominal fan power for $15$ November $2015$ (left) and $17$ November $2015$ (right).}
\label{fig:res_cap_setback}
\end{figure}

In this section, we present results relevant to the reserve scheduler. Two main factors that determine the amount of reserves are the building's energy capacity and the symmetry of reserve capacity. Apart from the physical properties of the building, the energy capacity depends also on the comfort zone's width. In this experiment, we specifically address the effect of enlarging the comfort zone during unoccupied hours to $12-35^\circ$C (the so-called night setback). We performed six full-day experiments with symmetric (equal up- and down-reserves) or asymmetric reserves, and with or without night setback. Note that the same price is assumed for up- and down-reserves in the asymmetric case.

Fig.~\ref{fig:res_cap_sym} shows the results for $20$ November when symmetric reserve capacities were assumed and for $21$ November when asymmetric capacities were used (in both days night setback was applied). The capacities are reported in $\%$ of the fan rated power ($2500$ W). The reserve capacity is maximized at night when the comfort zone is enlarged, and during the hottest part of daytime. In case of symmetric reserves, the maximum capacity is slightly less than $40\%$ of the rated fan power. In case of asymmetric reserves, the maximum up-reserve capacity is approximately $15\%$, whereas the down-reserve capacity is more than $60\%$. These experimental results are in agreement with relevant simulation results in \cite{VrettosIFAC2014,VrettosTPS2016}, and show that down-reserves (consumption increase) are preferable for commercial buildings equipped with energy-efficient controllers, because down-reserves can be provided without increasing the baseline consumption and the energy cost.

Fig.~\ref{fig:res_cap_setback} compares the experimental results of $15$ November when setback was used with those of $17$ November when no setback was applied (the reserve capacities were symmetric in both dates). With setback most reserve is provided at night, whereas without setback the reserve provision coincides with the highest cooling load in the middle of the day. Although the experiment was conducted with setback and symmetric reserves both on $15$ and $20$ November, the capacity profiles during daytime are considerably different due to different weather conditions and building models.

\begin{table}[t]
\renewcommand{\arraystretch}{1.1}
\caption{Experimental (bold) and simulated (normal font) daily average reserve capacities in $\%$ of fan nominal power ($2500$ W).}
\centering
\begin{tabular}{c|cc|cc|cc|cc}
\hline
Date & \multicolumn{2}{c|}{\begin{tabular}{@{}c@{}}Symmetric, \\ Setback\end{tabular}} & \multicolumn{2}{c|}{\begin{tabular}{@{}c@{}}Asymmetric, \\ Setback\end{tabular}} & \multicolumn{2}{c|}{\begin{tabular}{@{}c@{}}Symmetric, \\ No setback\end{tabular}} & \multicolumn{2}{c}{\begin{tabular}{@{}c@{}}Asymmetric, \\ No setback\end{tabular}}\\
& $R_\text{u}$ & $R_\text{d}$ & $R_\text{u}$ & $R_\text{d}$ & $R_\text{u}$ & $R_\text{d}$ & $R_\text{u}$ & $R_\text{d}$\\
\hline
$15/11$ &  $\mathbf{15.61}$ & $\mathbf{15.61}$ & $6.45$  & $26.15$ & $1.91$ & $1.91$ & $1.64$ & $4.42$ \\
\hline
$16/11$ & $9.09$ & $9.09$ &  $\mathbf{3.90}$ & $\mathbf{14.88}$  & $0.74$  & $0.74$  & $0.90$ & $1.85$ \\
\hline
$17/11$ & $11.44$  & $11.44$  & $5.39$ &  $21.21$ & $\mathbf{2.24}$ & $\mathbf{2.24}$ & $1.85$  & $4.59$ \\
\hline
$18/11$ & $16.70$ & $16.70$  & $7.78$  & $31.75$  & $3.94$  & $3.94$  & $\mathbf{3.28}$ & $\mathbf{9.91}$ \\
\hline
$20/11$ & $\mathbf{28.95}$ & $\mathbf{28.95}$ & $11.89$ & $49.66$ & $15.07$ & $15.07$ & $7.82$ & $28.81$ \\
\hline
$21/11$ & $22.10$  &  $22.10$  &  $\mathbf{10.72}$ & $\mathbf{46.55}$ & $13.60$ & $13.60$  & $6.79$ & $26.30$ \\
\hline
\end{tabular}
\label{tab:res_cap}
\end{table}

To have a fair comparison under the same external conditions, we simulated the reserve capacity scheduling for all combinations of symmetry and setback using the building model. The simulation and experimental results are shown in Table~\ref{tab:res_cap}. The capacity ranges from low values below $1\%$ to high values nearly $50\%$, and it heavily depends on reserve symmetry, setback, and weather conditions. The night setback increases the capacity by $177.0\%$ on average for symmetric reserves, by $107.0\%$ for asymmetric up-reserves, and by $150.7\%$ for asymmetric down-reserves. If setback is already used, adopting asymmetric capacities instead of symmetric capacities reduces the up-reserves by $55.6\%$ but increases the down-reserves by $83.1\%$, and so the net effect is an increase of $13.7\%$ in the total capacity.

\section{Room Climate Control (Level 2)} \label{results_lv2}
\subsection{Comfort Satisfaction}
\begin{figure}[t]
\centering \includegraphics[width=0.95\textwidth]{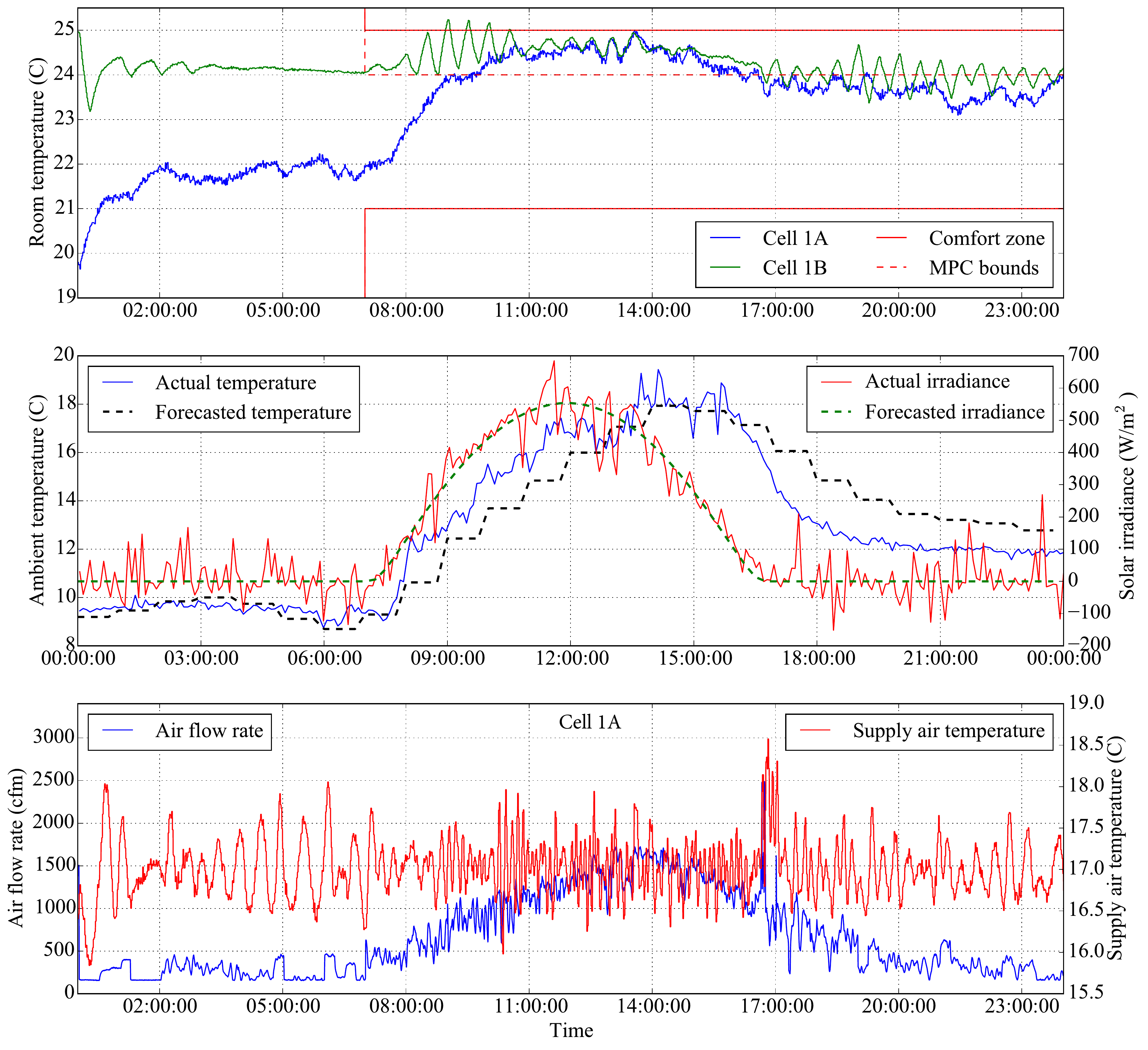}
\caption{Experimental results for the room climate controller under frequency regulation on $17$ November (symmetric reserves, without night setback).}
\label{fig:level2_comfort_res_17Nov}
\end{figure}

\begin{figure}[t]
\centering \includegraphics[width=0.95\textwidth]{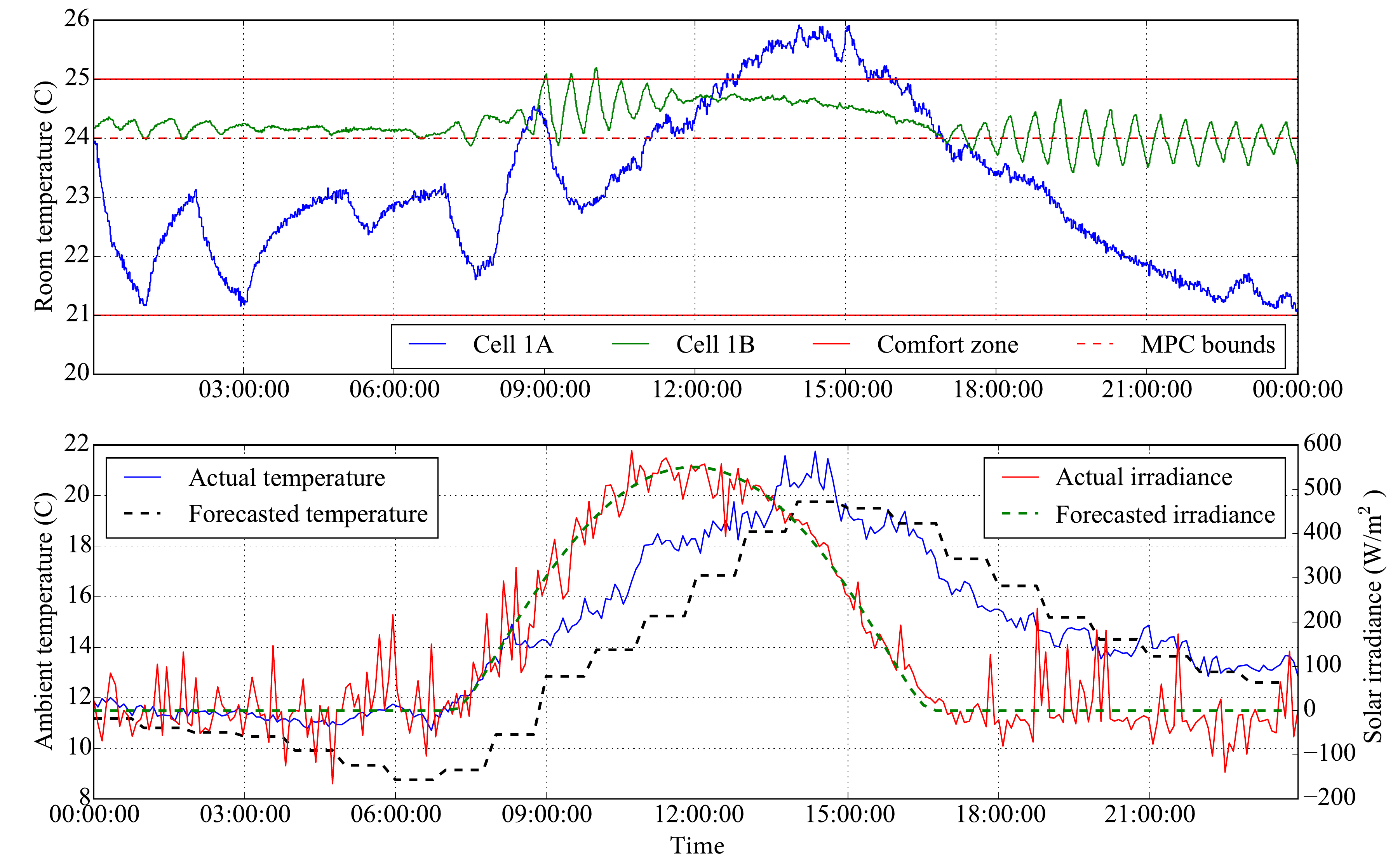}
\caption{Experimental results for the room climate controller under frequency regulation on $18$ November (asymmetric reserves, without night setback).}
\label{fig:level2_comfort_res_18Nov}
\end{figure}

\begin{figure}[t]
\centering \includegraphics[width=0.95\textwidth]{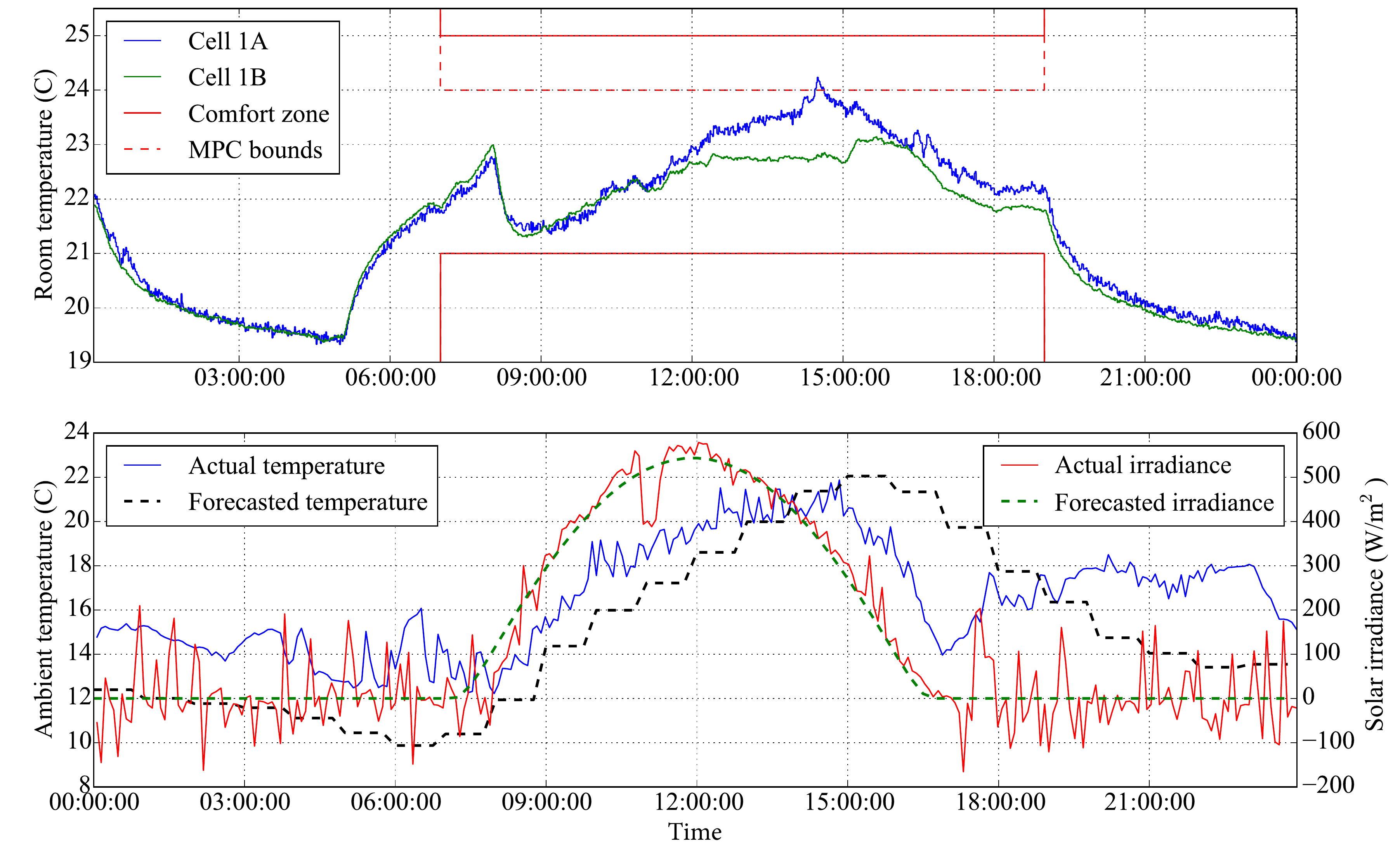}
\caption{Experimental results for the room climate controller under frequency regulation on $20$ November (symmetric reserves, with night setback).}
\label{fig:level2_comfort_res_20Nov}
\end{figure}

\begin{figure}[t]
\centering \includegraphics[width=0.95\textwidth]{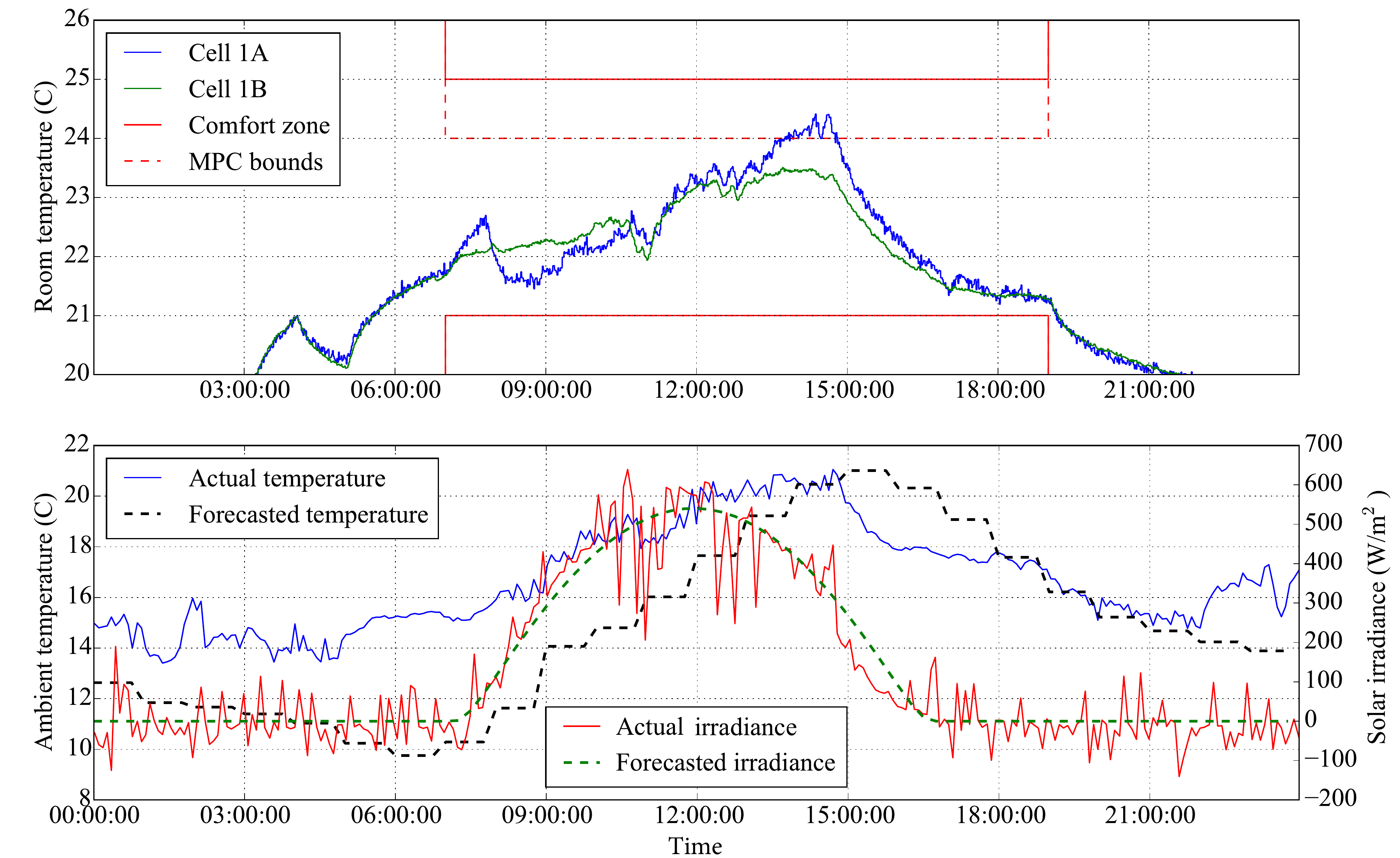}
\caption{Experimental results for the room climate controller under frequency regulation on $21$ November (asymmetric reserves, with night setback).}
\label{fig:level2_comfort_res_21Nov}
\end{figure}

Experimental results for $17$ November are shown in Fig.~\ref{fig:level2_comfort_res_17Nov}, where the top plot shows the temperature trajectories in cells 1A and 1B, the middle plot presents the forecasts and actual values for ambient temperature and solar irradiance\footnote{The total global irradiance is shown, which includes the long-wave radiation losses from the building envelope to the atmosphere, and it can be negative at night. This effect is known as nighttime radiation cooling \cite{givoni1994passive}.}, whereas the bottom plot shows the SAT and the air flow rate in cell 1A. The comfort zone is indicated with red: the actual upper limit (red solid line) is $25^\circ$C, but a tighter limit of $24^\circ$C (red dashed line) is used within the MPC to account for modeling and prediction errors. Similar results for $18$, $20$, and $21$ November are shown in Figs.~\ref{fig:level2_comfort_res_18Nov}, \ref{fig:level2_comfort_res_20Nov}, and \ref{fig:level2_comfort_res_21Nov}, but without including the SAT and air flow rate plot due to space limitations.

In Figs.~\ref{fig:level2_comfort_res_17Nov} and \ref{fig:level2_comfort_res_18Nov} the cell 1B is under energy efficient operation and the temperature remains close to the upper limit of the comfort zone. On the other hand, in Figs.~\ref{fig:level2_comfort_res_20Nov} and \ref{fig:level2_comfort_res_21Nov} the cell 1B is in a ``regulation-ready'' operation mode, namely the consumption of the \ac{HVAC} system is scheduled identically to cell 1A to allow reserve provision, but no regulation signal is received. For this reason, the temperature trajectories of the two cells are very close to each other for most of $20$ and $21$ November.\footnote{The discrepancies from $12.00$ to $19.00$ on $20$ November are due to the calibration differences between the fan models of the two cells (see Section~\ref{exp_settings}). The discrepancies from $07.00$ to $17.00$ on $21$ November are because of interruptions in the hierarchical control in cell 1B due to server connection timeout error from approximately $07.00$ to $11.00$. When the server was unresponsive, the cell was controlled by an existing fallback controller.}

The temperature trajectory of cell 1B in Figs.~\ref{fig:level2_comfort_res_17Nov} and \ref{fig:level2_comfort_res_18Nov} remains mostly in the band $[24-25]^\circ$C, which illustrates the necessity of tightening the comfort zone constraints in the MPC to compensate for modeling errors. The temperature trajectory of cell 1A is more variable and it follows the scheduled reserve and air flow rate. On $17$ November (Fig.~\ref{fig:level2_comfort_res_17Nov}) frequency regulation is provided while respecting the comfort zone.

\begin{figure}[t]
\centering
\begin{minipage}{0.49\linewidth}%
\centering \includegraphics[width=1\textwidth]{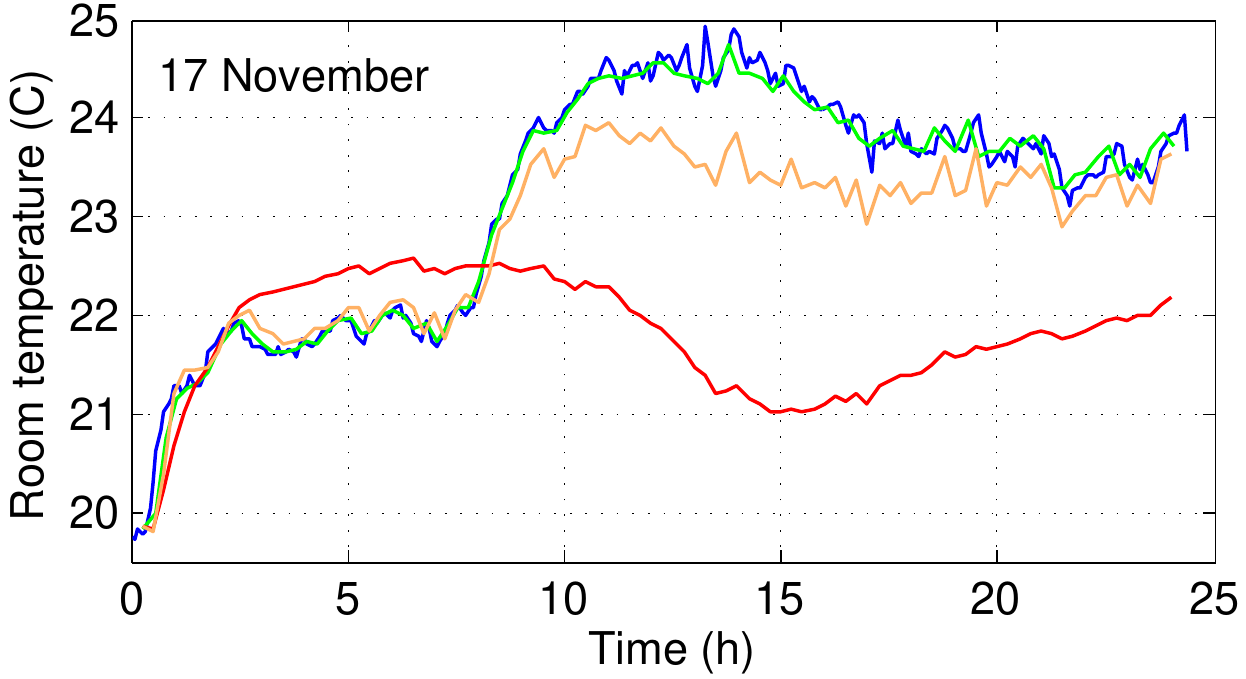}
\end{minipage}
\begin{minipage}{0.49\linewidth}
\centering \includegraphics[width=1\textwidth]{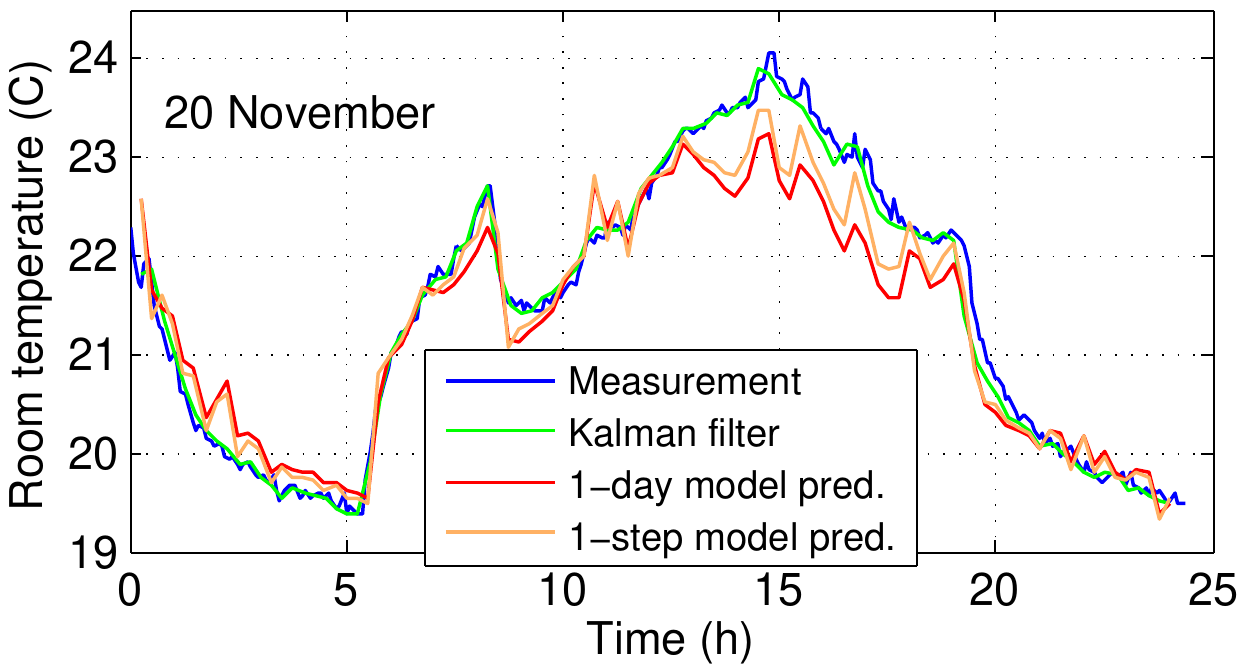}
\end{minipage}
\caption{Model and Kalman filter performance. Left: Results for $17$ November with the older model. Right: Results for $20$ November with the new model.}
\label{fig:model_and_KF_performance}
\end{figure}

\begin{figure}[h]
\centering \includegraphics[width=0.95\textwidth]{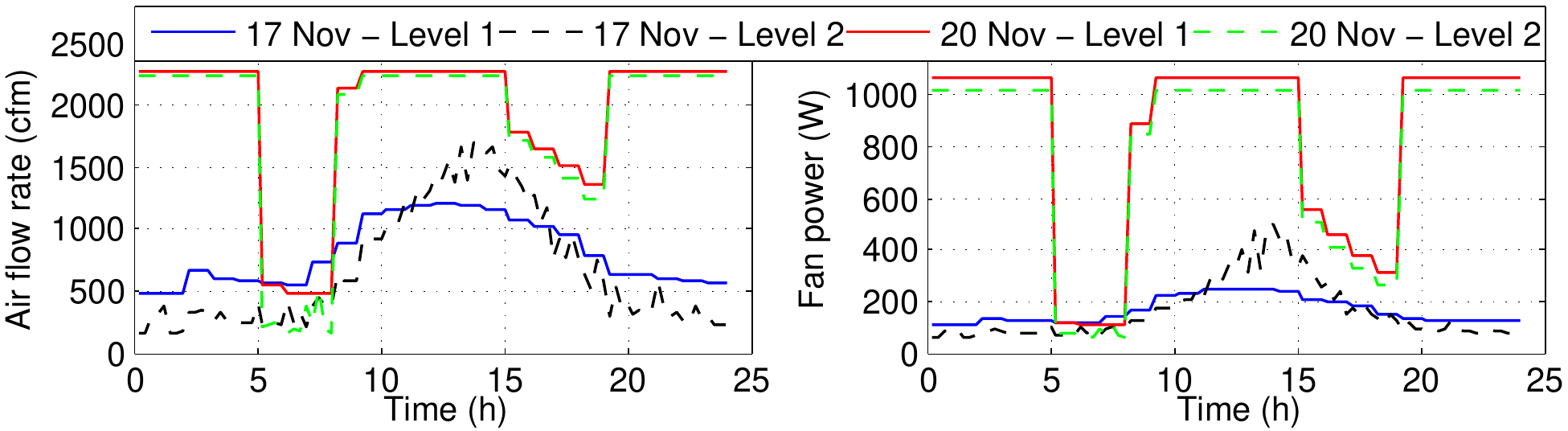}
\caption{Left: Air flow rate schedule in level $1$ and level $2$. Right: Fan power schedule in level $1$ and level $2$.}
\label{fig:mpc_reaction_model_error}
\end{figure}

However, on $18$ November (Fig.~\ref{fig:level2_comfort_res_18Nov}) the comfort zone is violated from $13.00$ to $16.00$ in cell 1A, but not in cell 1B. This happens because: (i) the ambient temperature is higher than the day-ahead forecast from the beginning of the day until $15.00$, and (ii) asymmetric reserves are used. The asymmetry allows for a more aggressive scheduling with a larger down-reserve capacity on $18$ November in comparison with $17$ November when symmetric reserves are used (see Table~\ref{tab:res_cap}).

The control performance is significantly better on $20$ and $21$ November (Figs.~\ref{fig:level2_comfort_res_20Nov} and \ref{fig:level2_comfort_res_21Nov}) despite the large discrepancies between the day-ahead ambient temperature forecasts and the actual values. The improvement is due to the recently calibrated building model (see Section~\ref{experiment_plan}). No comfort zone violations occur and moreover the temperature is below the MPC constraint of $24^\circ$C for most of the time. Therefore, periodic model calibration (for example on a weekly or daily basis) is important to account for seasonality and eliminate systematic errors.

These results show that if the model and weather predictions are sufficiently accurate, the robust reserve scheduler allows a commercial building to bid in day-ahead markets for frequency regulation. On the other hand, if the modeling and prediction errors exceed the controller's robustness margin, reserve provision for frequency regulation might have an adverse effect on occupant comfort.

The temperature trajectory in Figs.~\ref{fig:level2_comfort_res_20Nov} and \ref{fig:level2_comfort_res_21Nov} is typical for a building with a night setback. The controller chooses to overcool the space at night in order to generate higher revenue by offering a larger reserve capacity. In contrast, the reserve capacity is smaller during working hours and the room temperature is higher. A comparison of the temperature trajectories in cells 1A and 1B shows that tracking the RegD signal has little effect on room temperature due to the signal's limited energy content.

\subsection{Model and Estimator Performance}
Fig.~\ref{fig:model_and_KF_performance} compares the out-of-sample performance of the older model (left plot) and the new model (right plot). The blue curve is the measured room temperature, the green curve is the estimated temperature with the Kalman filter, the red curve corresponds to a day-ahead model prediction, whereas the orange curve shows the step-ahead temperature predictions. Clearly, the new model outperforms the older one, especially for the day-ahead predictions. This is why the performance of the level $2$ controller is much better on $20$ November than on $17$ November in terms of comfort zone violations.

The effect of model accuracy on MPC operation is shown in Fig.~\ref{fig:mpc_reaction_model_error}. On $17$ November the model mismatch is large, which results in a significant discrepancy in the scheduled air flow rate and fan power between levels $1$ and $2$. The MPC reacts on the modeling error by reducing the cooling power in level $2$ during night hours and increasing it during daytime. In this way, the MPC provides the same amount of electric reserve in daytime with less change in air flow rate by taking advantage of the nonlinear fan curve. On the other hand, the model mismatch is small on $20$ November and so the air flow and fan power schedules of level $1$ and level $2$ are similar. In fact, level $2$ consistently schedules less cooling power than level $1$ because the air flow constraints are relaxed \cite[Equations 15, 24]{VrettosTSG2016Exp_p1}, and the reserve scheduling in level $1$ is robust and thus conservative.

\subsection{Fan Heat Gain at High Speeds}
\begin{figure}[t]
\centering \includegraphics[width=0.95\textwidth]{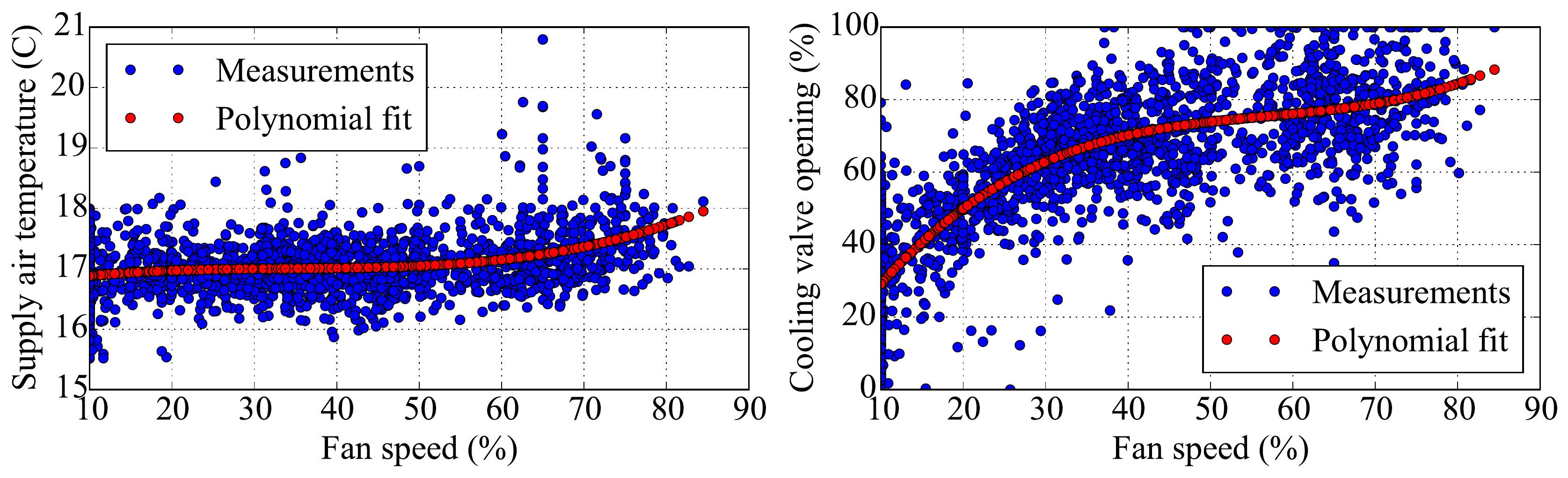}
\caption{The dependence of SAT and cooling valve opening on fan speed.}
\label{fig:fan_heat_gain}
\end{figure}

We present results on the dependence of SAT and cooling valve opening on fan speed in Fig.~\ref{fig:fan_heat_gain}, where the blue points  are measurements and the red trend is a polynomial fit on them. As expected, the trend in cooling valve opening is increasing because the higher the fan speed the more cooling is required from the chilled water loop. The trend in SAT is a flat line for fan speeds up to $50\%$. However, for speeds above $50\%$ (and especially above $70\%$) there is a clear increasing SAT trend despite the increased cooling valve opening.

These results lead to an interesting observation: the heat gain due to the rotation of the fan is significant at high speeds and it cannot be effectively rejected by exchanging heat with the chilled water loop. According to Fig.~\ref{fig:fan_heat_gain}, if the fan operates at a $70\%$ speed or higher, the SAT will likely have a steady-state deviation from the setpoint $17^\circ$C that can be as high as $1^\circ$C. Steady-state SAT deviations might result in comfort zone violations, because the controller assumes the SAT fixed to $17^\circ$C. This did not create problems in our experiment because the scheduled fan speed by the MPC was at most $70\%$.

\subsection{Effect on Energy Consumption}
A major concern when providing \acp{AS} with commercial buildings is the effect on energy consumption. Ref.~\cite{Backhaus2015} reported a round-trip efficiency of $46\%$ when a building responded to demand response events in an experiment. There are two types of efficiency losses relevant to frequency regulation: ``reserve availability efficiency loss'' and ``reserve utilization efficiency loss'' \cite{VrettosTSTE2016}. The first one is the efficiency loss due to scheduling the consumption in an energy suboptimal way to be able to provide frequency reserves, if requested. The second one is the additional efficiency loss while tracking the frequency regulation signal.

We report efficiency results in Table~\ref{tab:efficiency} for: (i) $15$-$18$ November when the cell 1B was under energy efficient operation (to quantify the reserve availability efficiency loss); and (ii) $20$ November when cell 1B was in regulation-ready operation mode (to quantify the reserve utilization efficiency loss). The efficiency loss is calculated comparing the energy consumption of cell 1A with that of the benchmark cell 1B. We use two different definitions of energy consumption: (i) electric energy consumption of the fan, and (ii) thermal cooling power consumption of each cell. The latter is calculated based on the chilled water flow rate ($\dot{m}_\text{cw}$), as well as the supply ($T_\text{ch,s}$) and return ($T_\text{ch,r}$) chilled water temperatures using
\begin{align} \label{Pcool_def}
    P_\text{cool} = \dot{m}_\text{cw} \cdot \left(T_\text{ch,r}-T_\text{ch,s}\right)~.
\end{align}

\begin{table}[t]
\renewcommand{\arraystretch}{1.1}
\caption{Effect of frequency regulation on energy consumption}
\centering
\begin{tabular}{c||cc|cc|cc}
\hline
\begin{tabular}{@{}c@{}}1B operation\\ mode\end{tabular} & \multicolumn{2}{c|}{\begin{tabular}{@{}c@{}}Energy efficient\\ $15/11$-$18/11$\end{tabular}}   &  \multicolumn{4}{c}{\begin{tabular}{@{}c@{}}Regulation-ready\\ $20/11$ ($0$-$24$~h) $20/11$ ($0$-$12$~h)\end{tabular}} \\
 \hline
Cell & 1A & 1B & 1A  & 1B & 1A  & 1B \\
 \hline
Fan energy (kWh) & $27.22$  & $16.23$ & $19.24$ & $20.91$ & $9.40$ & $10.19$ \\
Cooling (gpm $\cdot$ F) & $1989.92$ & $1800.55$ & $772.16$ & $788.69$ & $328.78$ & $331.09$\\
Mean temp. ($^\circ$C) & $22.85$ & $24.43$ & $21.45$ & $21.26$ & $21.05$ & $21.05$\\
\hline
\end{tabular}
\label{tab:efficiency}
\end{table}

Based on the results of Table~\ref{tab:efficiency}, the reserve availability efficiency loss is significant and equal to approximately $68\%$ in terms of fan consumption and $11\%$ in terms of cooling power from the chiller. However, the additional consumption in cell 1A is not entirely wasted because it results in a lower average temperature. When cell 1B is in the regulation-ready mode, the cell 1A consumes less energy than cell 1B despite frequency regulation. The non-negligible difference in the average temperature of the two cells is due to the imperfections in fan model calibration and the limited temperature sensor accuracy. However, even from $00.00$ to $12.00$ when both average cell temperatures are $21.05^\circ$C, the consumption of cell 1A is still lower than that of cell 1B. This result indicates that the ``reserve utilization efficiency loss'' is negligible while tracking a fast-moving regulation signal like RegD.

\subsection{MPC Computation Time}
\begin{figure}[t]
\centering \includegraphics[width=0.85\textwidth]{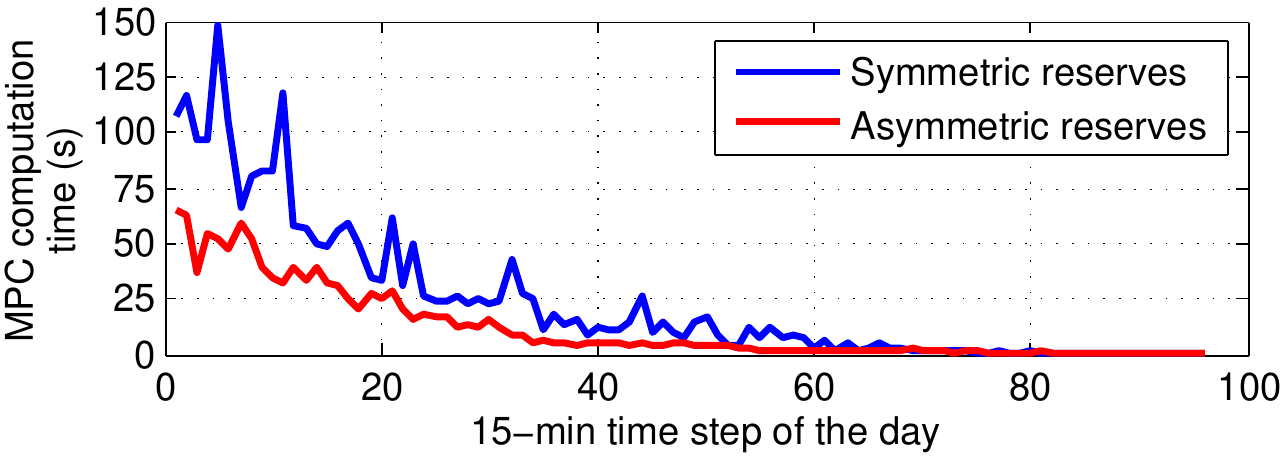}
\caption{The average MPC computation time depending on the time of the day and on reserve symmetry.}
\label{fig:mpc_comp_time}
\end{figure}

The MPC computation time is sufficiently low for our demonstration. As shown in Fig.~\ref{fig:mpc_comp_time}, the longest computation time is $150$~s for symmetric reserves and $65$~s for asymmetric reserves. The computation time for asymmetric reserves is lower because the problem is simpler and smaller \cite[Section IV.D]{VrettosTSG2016Exp_p1}.

The computation time decreases at the end of the day because a reducing MPC horizon is used. After the $70^\text{th}$ time step, when the MPC prediction horizon is smaller than $26$ time steps ($6.5$ hours), the computation time is less than $2$~s. Therefore, Fig.~\ref{fig:mpc_comp_time} can be used to select the prediction horizon's length depending on the maximum allowable computation time. Since the computation time grows exponentially with the number of variables of the nonlinear optimization problem, a shorter prediction horizon might be necessary for larger buildings.

\section{Regulation Signal Tracking (Level 3)} \label{results_lv3}
\subsection{Control Performance Metrics}
In this section, we present results from level $3$ and evaluate the tracking performance of the regulation signal. The following metrics are used

\begin{align}
e_{\text{t},k} &= e_{\text{c},k}/P_{\text{d},k},~~ e_{\text{c},k} = P_{\text{d},k}-P_{\text{f},k} \label{track_perc_error}\\
e_{\text{r},k} &= \begin{cases}
                    e_{\text{c},k}/R_{\text{u},k}, & \textrm{if} \hspace{2mm} w_k<0\\
                    e_{\text{c},k}/R_{\text{d},k}, & \textrm{if} \hspace{2mm} w_k\geq0
                  \end{cases} \label{res_perc_error}\\
e_\text{me} &= (1/N_\text{exp}) \cdot \sum\nolimits_{k=0}^{N_\text{exp}-1} e_{\text{c},k}\\
e_\text{mae} &= (1/N_\text{exp}) \cdot \sum\nolimits_{k=0}^{N_\text{exp}-1} \left|e_{\text{c},k}\right|\\
e_\text{rmse} &= \sqrt{(1/N_\text{exp}) \cdot \sum\nolimits_{k=0}^{N_\text{exp}-1}e_{\text{c},k}^2}\\
e_\text{t,mape} &= (1/N_\text{exp}) \cdot \sum\nolimits_{k=0}^{N_\text{exp}-1} \left|e_{\text{t},k}\right| \\
e_\text{r,mape} &= (1/N_\text{exp}) \cdot \sum\nolimits_{k=0}^{N_\text{exp}-1} \left|e_{\text{r},k}\right|~,\label{res_mape}
\end{align}
where $P_{\text{f},k}$ denotes the instantaneous fan power, $w_k$ denotes the normalized regulation signal, and $N_\text{exp}$ denotes the experiment duration. The metrics $e_{\text{t},k}$ and $e_{\text{r},k}$ are relative instantaneous errors but the normalization is performed using the desired fan power $P_{\text{d},k}$ in $e_{\text{t},k}$, and the up- ($R_{\text{u},k}$) or down-reserve capacity ($R_{\text{d},k}$) in $e_{\text{r},k}$. The mean error $e_\text{me}$ is used to measure any biases in the control response, whereas $e_\text{mae}$ is the \ac{MAE} during the experiment. The \ac{RMSE} $e_\text{rmse}$ penalizes more the large control errors, for example due to overshoots and undershoots. The metric $e_\text{t,mape}$ is the tracking \ac{MAPE}, and $e_\text{r,mape}$ is the reserve \ac{MAPE}. We use the metric $e_\text{r,mape}$ because it describes the relative size of the control error with respect to the reserve capacity.

In addition, we use the score proposed by PJM for evaluating the performance of frequency regulation. The total score $S_\text{tot}$ consists of three parts, namely the correlation score $S_\text{c}$, the delay score $S_\text{d}$ and the precision score $S_\text{p}$, which are defined as \cite{PJMperformanceScores}
\begin{align}
S_\text{c} &= \underset{\tau \in [0, 5~\text{min}]}{\max} (R_\text{cor}) \\
S_\text{d} &= \Big|\frac{\tau^* - 5~\text{min}}{5~\text{min}}\Big|,~\tau^* = \underset{\tau \in [0, 5~\text{min}]}{\text{argmax}} (R_\text{cor}) \\
S_\text{p} &= 1-(1/n) \cdot \sum\nolimits_{k=1}^{n} \Big|\frac{e_{\text{c},k}}{\bar{P}_{\text{d},\text{h}}}\Big| \\
S_\text{tot} &= (1/3) \cdot S_\text{c} + (1/3) \cdot S_\text{d} + (1/3) \cdot S_\text{p}~.
\end{align}

The correlation score is the maximum correlation $R_\text{cor}$ of the desired power $P_{\text{d},k}$ and fan power $P_{\text{f},k}$, and $\tau^*$ is the time shift at which the correlation is maximized ($\tau$ takes a value from $0$ to $5$~min with a step of $10$ s). We calculate the delay score based on the time shift with maximum correlation. In the precision score calculation, we normalize the absolute control error by the average hourly value of the reference signal $\bar{P}_{\text{d},\text{h}}$, whereas the total score is a weighted sum of the individual scores.

\subsection{Experimental Time Series Results}
\begin{figure}[t]
\centering \includegraphics[width=0.68\textwidth]{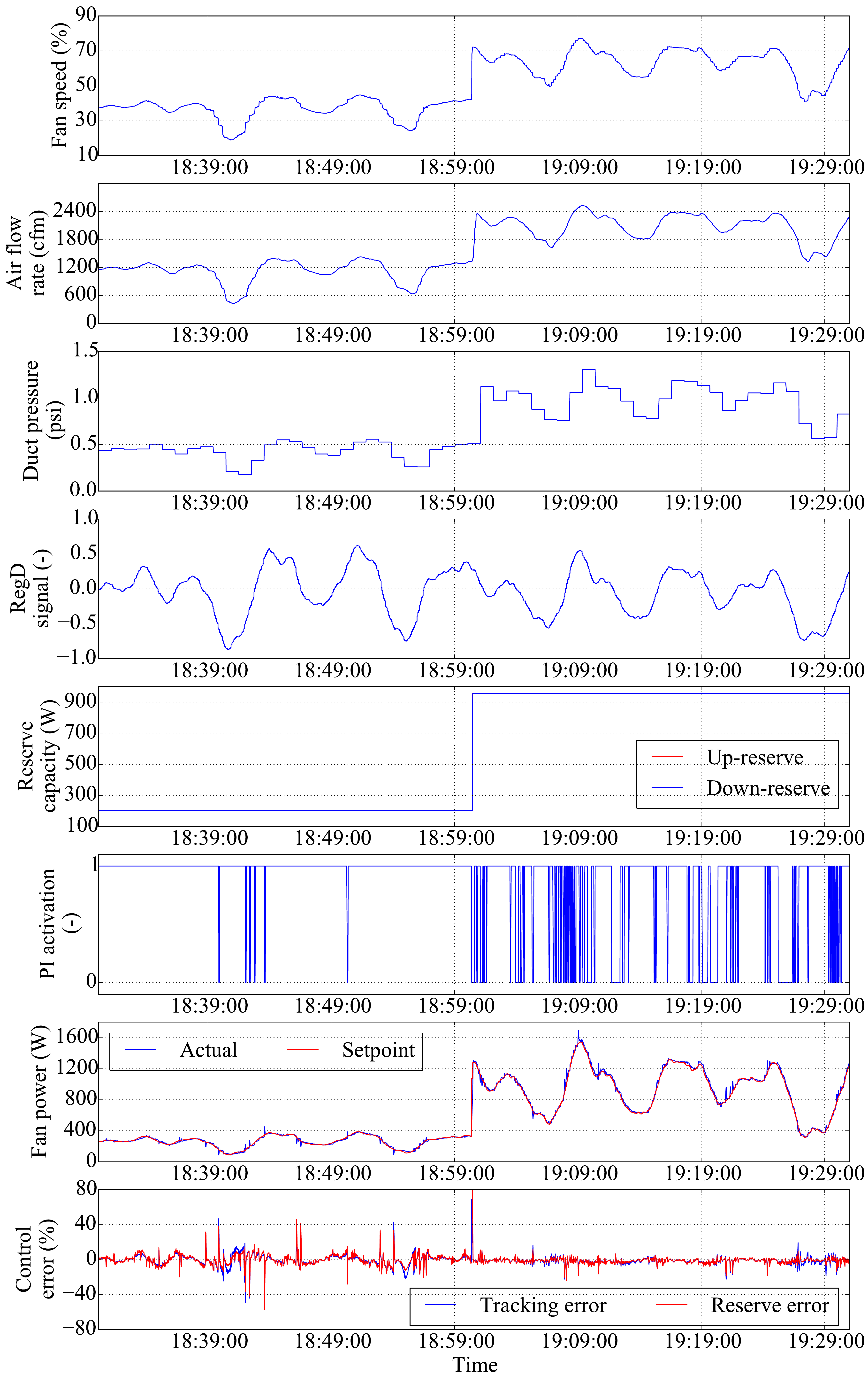}
\caption{Fan control and RegD signal tracking for a period of $1$ hour.}
\label{fig:ctrlFanLevel3}
\end{figure}

In Fig.~\ref{fig:ctrlFanLevel3} we present results from the operation of level $3$ controller from 18.30 to 19:30 on $20$ November $2015$. The duct pressure is quadratic to the fan speed, as expected from the fan laws. Since the duct system is designed to sustain the pressure corresponding to the maximum fan speed, and because the fan speed does not exceed its maximum value ($90\%$) during frequency regulation, pressure constraints were not necessary in the reserve scheduling and MPC formulations in our experiment.

The RegD signal changes direction very often and has a limited energy content. During periods of time when the RegD signal is relatively flat, or the reserve capacity is low, the \ac{PI} controller is active. On the other hand, whenever the changes in fan power are rapid, the control switches to the model-based feedforward control.

The tracking of the RegD signal is generally very good. However, when large rapid changes in fan power are requested, overshoots or undershoots might appear. In addition, if the reserve capacities change significantly at the beginning of each full hour, temporarily large errors might occur. In general, the instantaneous percentage errors $e_{\text{t},k}$ and $e_{\text{r},k}$ are higher at a low operating fan power and low reserve capacity.

\subsection{Evaluation of Tracking Performance}
The performance metrics \eqref{track_perc_error}-\eqref{res_mape} for the $6$ days of Table~\ref{tab:res_cap} are presented in Table~\ref{tab:tracking_error}. $e_\text{r,mape}$ is larger than $e_\text{t,mape}$ because small reserve capacities are offered for a large part of the experiment. The mean error $e_{\text{me}}$ has a negative bias, which means that the fan power is more often higher than the desired setpoint because the control overshoots are larger than the undershoots.

\begin{figure}[t]
\centering \includegraphics[width=0.95\textwidth]{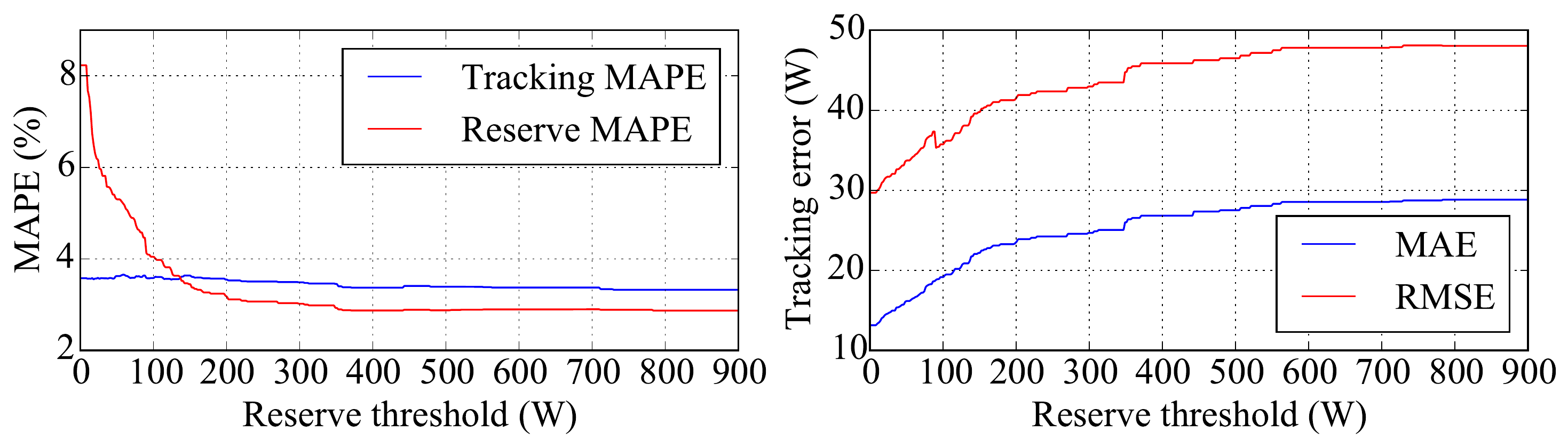}
\caption{Dependence of tracking and reserve errors on the reserve threshold.}
\label{fig:resErrVsResThres}
\end{figure}

We investigate the dependence of control performance on the minimum reserve capacity, which we call ``reserve threshold'' and denote by $R_\text{thr}$. The metrics $e_\text{mae}$, $e_\text{rmse}$, $e_\text{t,mape}$ and $e_\text{r,mape}$ are recalculated considering only the time steps when $R_{\text{u},k}\geq R_\text{thr}$ if $w_k<0$, and $R_{\text{d},k}\geq R_\text{thr}$ if $w_k\geq0$. We repeat this procedure for different $R_\text{thr}$ values and present the results in Fig.~\ref{fig:resErrVsResThres}.

In contrast to $e_\text{t,mape}$, $e_\text{r,mape}$ decreases rapidly as $R_\text{thr}$ increases in the range $[0, 200]$~W. This happens because (for the same absolute control error) $e_{\text{r},k}$ decreases if $R_{\text{u},k}$ or $R_{\text{d},k}$ increase. On the other hand, $e_\text{mae}$ and $e_\text{rmse}$ generally increase as $R_\text{thr}$ increases because the higher the reserve capacity the larger the fan power change, and thus the higher the errors due to the overshoots and undershoots. Fig.~\ref{fig:resErrVsResThres} can provide us with bounds on the reserve capacity from a tracking performance point of view.

In Table~\ref{tab:pjm_scores} we report the PJM scores calculated separately for the period $15$ - $18$ November and the period $20$ - $21$ November. Different scores are calculated for each hour (only if the reserve capacity is non-zero) \cite{PJMperformanceScores}, and the values in Table~\ref{tab:pjm_scores} are hourly averages. The frequency regulation performance is exceptional during the whole experiment. For comparison, the highest possible total score is $S_\text{tot}=1$ and the minimum $S_\text{tot}$ accepted by PJM is $0.75$. The scores are slightly higher for $20-21$ November because the building provides a larger reserve capacity compared with $15-18$ November.

\begin{table}[t]
\renewcommand{\arraystretch}{1.1}
\caption{RegD tracking performance metrics during the experiment}
\centering
\begin{tabular}{c|ccccc}
\hline
Metric & $e_\text{me}$ & $e_\text{mae}$ & $e_\text{rmse}$ & $e_\text{t,mape}$ & $e_\text{r,mape}$\\
Value & $-5.66$~W & $12.45$~W & $27.00$~W & $3.58\%$ & $8.23$\%\\
\hline
\end{tabular}
\label{tab:tracking_error}
\end{table}

\begin{table}[t]
\renewcommand{\arraystretch}{1.1}
\caption{PJM scores for tracking the RegD signal}
\centering
\begin{tabular}{c|cccc}
\hline
Score & $S_\text{c}$ & $S_\text{d}$ & $S_\text{p}$ & $S_\text{tot}$ \\
$15-18$ Nov. & $0.89$ & $0.97$ & $0.96$ & $0.94$ \\
$20-21$ Nov. & $0.96$ & $0.99$ & $0.98$  & $0.98$\\
\hline
\end{tabular}
\label{tab:pjm_scores}
\end{table}

\subsection{Effect on Supply Air Temperature}
\begin{figure}[t]
\centering \includegraphics[width=0.95\textwidth]{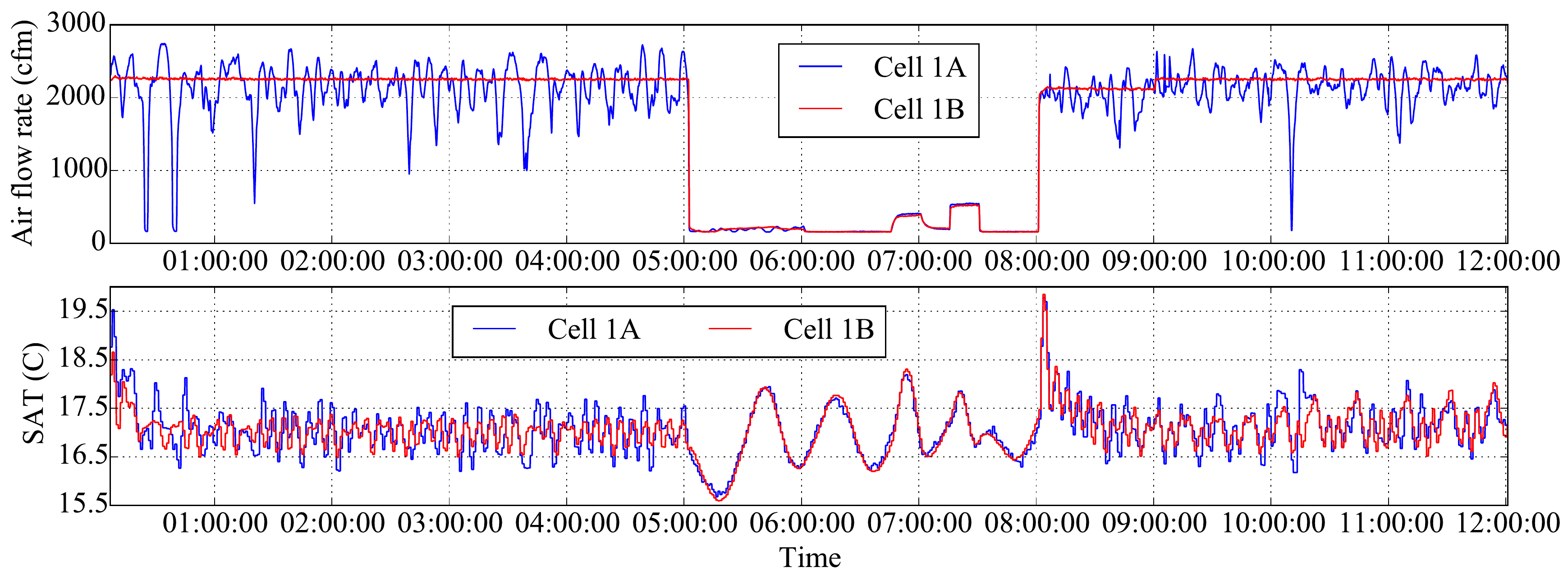}
\caption{Air flow rate and SAT in cells 1A and 1B on $20$ November.}
\label{fig:sat_fluctuation}
\end{figure}

Tracking the fast-moving RegD signal introduces high frequency oscillations on \ac{SAT} as shown in Fig.~\ref{fig:sat_fluctuation}. The SAT of cell 1A (frequency regulation) oscillates more than that of 1B (regulation-ready mode), especially after sudden changes in the regulation signal that induce sudden changes in the air flow rate. In addition, large excursions in SAT occur in both cells when the MPC changes the air flow setpoint significantly, for example at hour $08.00$. Moreover, the magnitude of SAT oscillations is high at low air flow rates, for example from $05.00$ to $08.00$.

\subsection{Effect of Fan Control on Chiller Power}
\begin{figure}[t]
\centering \includegraphics[width=0.95\textwidth]{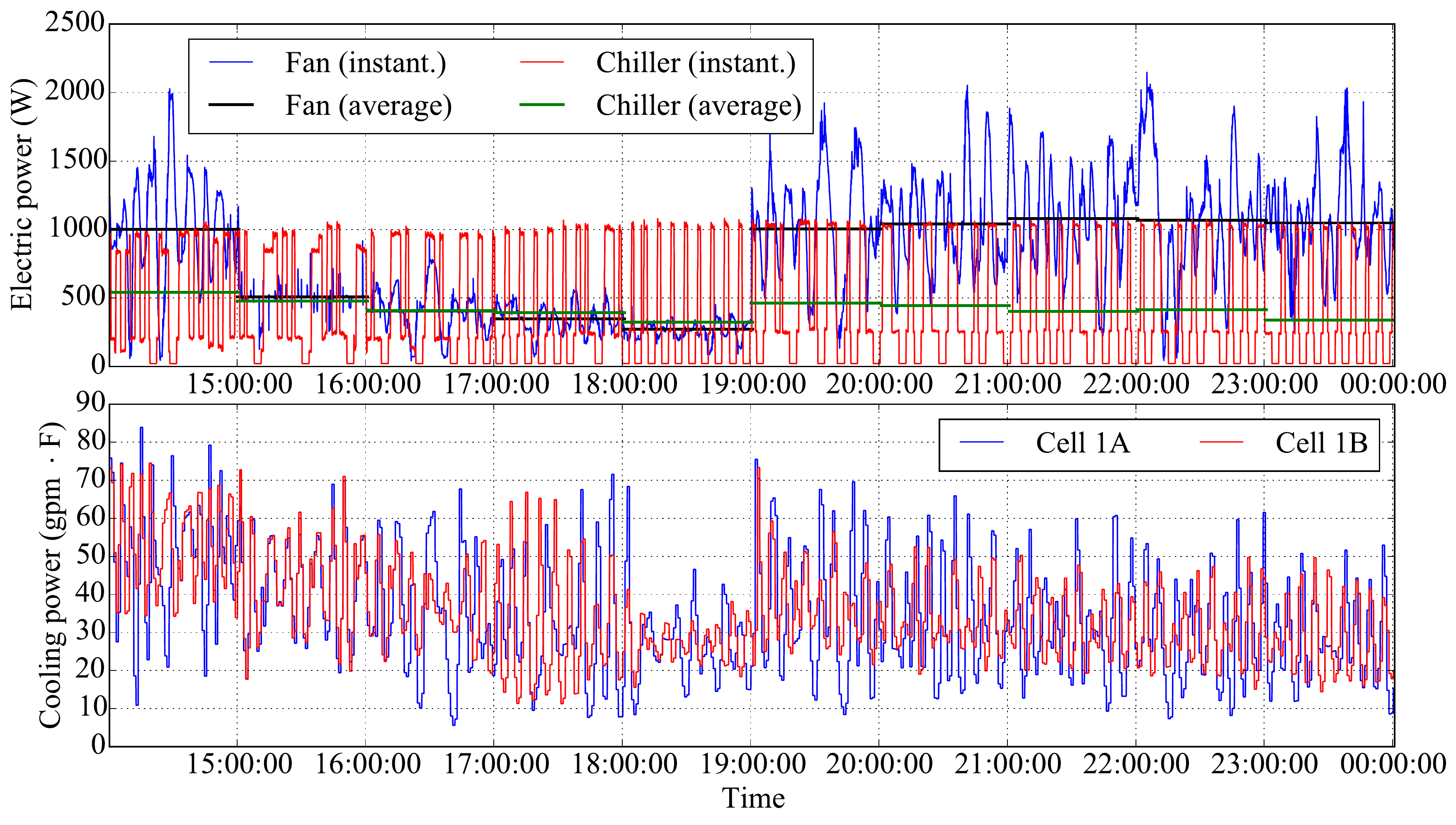}
\vspace{-0.2cm}
\caption{Effect of frequency regulation on the chiller and cooling power.}
\label{fig:chiller_effect}
\end{figure}

The fan and the chiller are thermally coupled through the chilled water loop, hence, it is worth investigating if the operation of the chiller is affected while providing frequency regulation with the fan. In Fig.~\ref{fig:chiller_effect} we present relevant experimental results for a duration of $10$ hours. The top plot shows the instantaneous and hourly-average electric power of the fan in cell 1A and the chiller. The bottom plot shows the cooling power in the chilled water loop for cells 1A and 1B calculated with \eqref{Pcool_def}.

The chiller has two stages and the electric power consumption is relatively constant at each stage. The chiller's cycling depends on the cooling load, which in turn depends on the fan power and ambient conditions. In general, as the fan power increases the chiller cycles more often and remains longer at the on state. This is shown in Fig.~\ref{fig:chiller_effect} where the average chiller power (green line) generally follows the average fan power (black line).

The effect of regulation is visible on the cooling power that fluctuates more in cell 1A compared with cell 1B (regulation-ready mode). Whenever the fan power increases, the cooling load also increases and the SAT tends to decrease. This is sensed by the SAT controller that opens the cooling valve to compensate for the SAT decrease, which in turn increases the cooling power in the chilled water loop. The delay in cooling power response depends on the time constant of the cooling valve's controller.

Despite the oscillations in cooling power, there is no observable effect on chiller's cycling and electric power. This happens because: (i) the chilled water is stored in a tank that provides some inertia; and (ii) the RegD signal is approximately zero-mean. Note that the gradual reduction in the hourly-average chiller electric power from 19.00 to 00.00 in Fig.~\ref{fig:chiller_effect} is mainly the result of a lower cooling need due to ambient temperature drop, rather than a side-effect of frequency regulation.

These results indicate that frequency regulation can be provided with fan control without side-effects on chiller consumption. However, this does not necessarily hold for regulation signals with a larger energy content such as RegA. In addition, chillers with continuous rather than duty-cycle control will likely display a more observable impact on consumption while providing frequency regulation, especially if there is no chilled water storage tank. In these cases, the level $3$ controller should be revised, which is an interesting direction for future work.

\begin{figure}[t]
\centering \includegraphics[width=0.85\textwidth]{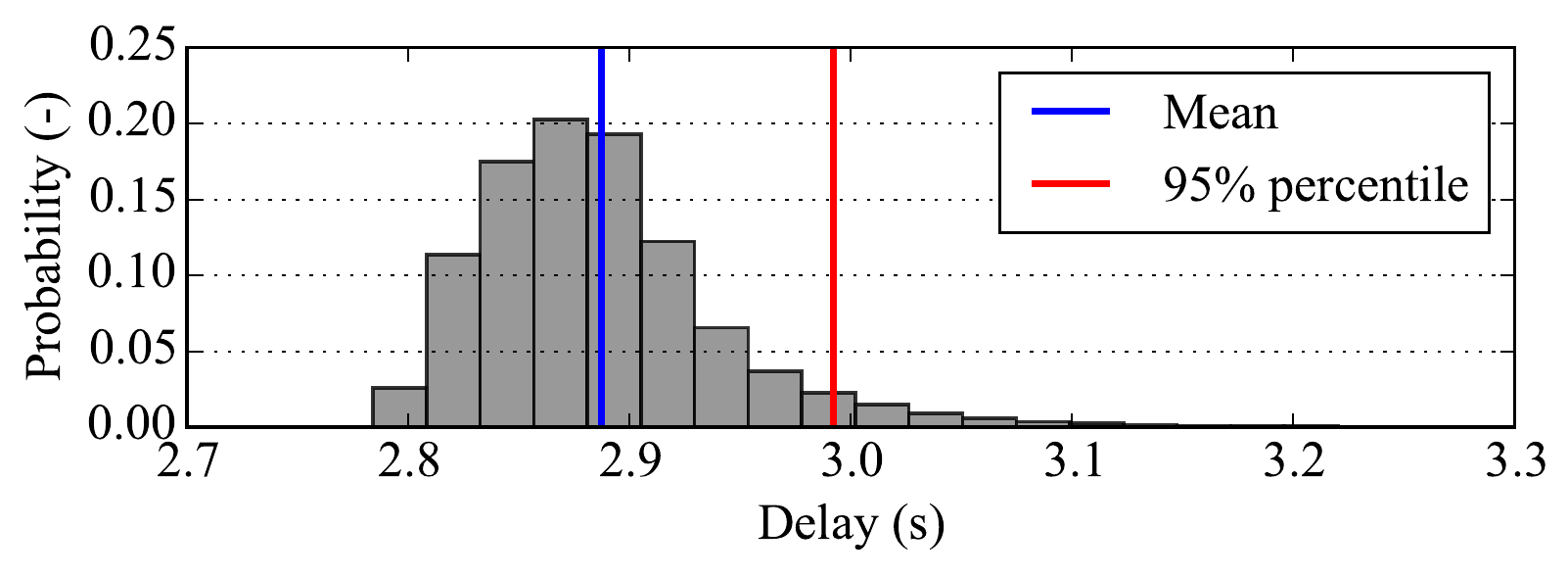}
\vspace{-0.2cm}
\caption{Histogram of communication delays during the experiment.}
\label{fig:delay_hist}
\end{figure}

\subsection{Analysis of Communication Delays}
A challenge in this experiment was the communication delays in measurements and actuation, which result in the overshoots and undershoots in fan power in Fig.~\ref{fig:ctrlFanLevel3}. In Fig.~\ref{fig:delay_hist} we present a histogram of the experienced delays during the whole experiment. The probability distribution of the delay is positively skewed with a mean value of $2.89$~s and a $95\%$-percentile of $2.99$~s. In fact, there exist a few very large delays in excess of $5$~s due to temporary unresponsiveness of the \ac{CWS}, which are not included in Fig.~\ref{fig:delay_hist}. Despite the fact that the average delay is large compared with the time step of level $3$ controller ($4$~s), the tracking performance of RegD signal is very good.

\section{Lessons Learned and Outlook} \label{sec_learnings}
\subsection{Lessons Learned}
\emph{Hierarchical control} is an efficient way to provide frequency regulation with commercial buildings because time-separated tasks are considered individually. Three control layers are essential: (i) a reserve capacity scheduler, (ii) a building climate controller to satisfy comfort while leaving enough slack for reserves, and (iii) a controller to track the regulation signal.

\emph{Frequency regulation accuracy:} High-quality frequency regulation can be provided by fan speed control. The RegD signal tracking is excellent even with large communication delays in the building automation system. A switched controller comprised of a feedforward controller and a \ac{PI} feedback controller with gain scheduling provides a fast response without compromising stability. This results in a total PJM score as high as $0.98$, which is well above PJM's limit of $0.75$.

\emph{Means to increase reserve capacity:} In our experiment, the fan provided $0.74-49.66\%$ of its rated power as reserve capacity, depending on ambient conditions and reserve assumptions. Allowing asymmetric reserve capacities and using a night setback are effective ways to increase the reserve potential from commercial buildings. In fact, down-reserves are preferable for buildings because the capacity can be offered without increasing baseline energy consumption.

\emph{Occupant comfort:} If the building bids in day-ahead \ac{AS} markets, respecting occupant comfort might be challenging if the building model and weather forecasts are not very accurate. Furthermore, asymmetric reserves result in a more aggressive scheduling that might increase comfort zone violations.

\emph{Building model:} An accurate building thermal model is essential for comfort satisfaction, especially in day-ahead AS markets. Periodic calibration of the building model helps to account for seasonality and eliminate offsets in modeling error.

\emph{Advantages of MPC:} Perhaps the most important advantage of MPC is that it identifies the optimal balance between reserve provision and energy efficiency. MPC additionally provides us with a baseline consumption ahead of real-time operation, which is beneficial from a practical point of view. Moreover, due to its predictive closed-loop nature it reacts to modeling and weather prediction errors in a way that minimizes occupant discomfort.

\emph{Robustness measures:} It is important to consider the regulation signal uncertainty when scheduling the reserve capacity. A conservative modeling of this uncertainty builds robustness to weather prediction and building modeling errors. Additional robustness can be obtained by tightening the comfort zone constraints in the MPC, and allowing a larger fan speed control band in the MPC compared with the reserve scheduler.

\emph{Effects of frequency regulation on building control:} Frequency regulation might introduce oscillations in SAT that can be reduced by appropriately tuning the cooling valve controller. In addition, if the MPC schedules the fan speed at very high values, the cooling loop might not be able to reject the additional heat gain due to fan rotation. On the positive side, there is little effect on the average energy consumption of the chiller while tracking an energy-constrained frequency regulation signal by controlling the fan power. However, the impact of fan control on chiller cycling may prevent the building from accurately following the regulation signal when measured against a baseline that includes the combined consumption of the fan and the chiller. This is an interesting area of further study.

\emph{Energy consumption:} Provision of frequency reserves entails some efficiency loss. The efficiency loss due to scheduling the HVAC consumption in a suboptimal way compared with an energy efficient building control can be as high as $67\%$. On the other hand, the efficiency loss while tracking frequency regulation signals with limited energy content is negligible.

\subsection{Outlook}
There are several avenues for follow-up work. Two direct extensions are to repeat the experiment with the RegA signal of PJM, which is slower but has more energy content, and/or with the heating loop of the \ac{AHU} enabled. In addition, performing the frequency regulation experiment using all four buildings of FLEXLAB will leverage the full potential of hierarchical control and verify the scalability of the approach.

In some HVAC systems a duct pressure controller regulates the pressure to a fixed setpoint. The combined operation of this controller and the dampers of each zone might reject the frequency regulation action \cite{zhao2013evaluation}. This is an important challenge that could not be addressed in this experiment at FLEXLAB as it requires testing on a large building.

The reserve scheduling optimization problem might be computationally heavy for buildings with many zones. An alternative is to approximate the nonlinear fan power curve with a piecewise affine function by introducing binary variables. The bilinear building dynamics can be approximated with sequential convex optimization \cite{Oldewurtel2011}, but the convergence is not guaranteed. Finally, the conservatism of the reserve scheduling problem can be reduced by generating scenarios from historical frequency regulation signals at the cost of reducing robustness.

\section{Conclusion} \label{sec_conclusion}
In Part II of this two-part paper, we reported experimental results for frequency regulation from a commercial building test facility (FLEXLAB). The results are very encouraging: the test building can track fast-moving signals such as RegD reliably, with very high accuracy, and with minimal effect on occupant comfort and the operation of the HVAC system. The results also indicate that a hierarchical control approach is appropriate for frequency regulation with day-ahead bidding of the reserve capacity, and it can be used in field tests and real-world implementations in larger buildings.

\bibliographystyle{elsarticle-num}
\bibliography{biblio_EV}

\end{document}